\documentclass[14]{article}
\usepackage[margin=1in]{geometry}
\usepackage{amsfonts,amsmath,amssymb}
\numberwithin{equation}{section}

\title{Fuzzy Limits of Fuzzy Functions}
\author{Abdulhameed Qahtan Abbood Altai \\ \small
University of babylon, Babil, Iraq 51002, \\ \small ahbabil1983@gmail.com}

\begin{document}

\maketitle

\paragraph{abstract}
In this paper, we propose the theory of fuzzy limit of fuzzy function depending on the  Altai's principle and using the representation theorem (resolution principle) to run the fuzzy arithmetic.

\paragraph{Keywords:}fuzzy limit of fuzzy function, two-sided fuzzy limits, one-sided fuzzy limits, fuzzy limit at infinite.

\section{Introduction and Preliminaries.} 
Zadeh introduced the concept of fuzzy set to assign to each object encountered in the real physical world that do not have precisely defined criteria of membership a grade of membership ranging between zero and one in 1965 [15]. Kramosil and Mich$\acute{a}$lek defined the concept of fuzzy metric space using continuous t-norms in 1975 [9]. The fuzzy metric spaces have very important applications in quantum physics, particularly, in connections with both string and $\epsilon^{(\infty)}$ theory which were studied by EI Naschie [14]. Matloka considered bounded and convergent sequences of fuzzy numbers and studied their properties in 1986 [11]. Sequences of fuzzy numbers also were discussed by Nanda $[13]$, Kwon $[10]$, Esi $[5]$ and many others. Burgin introduced the theory of fuzzy limits of functions based on the theory of fuzzy limits of sequences in 2000. He studied and developed the construction of fuzzy limits of functions similar to the one of the fuzzy limits of sequences based on the concept of $r-$limit of function $f$ [3]. In 2010, Altai defined the fuzzy metric spaces in a new way, that every real number $r \in \mathbb{R}$ is replaced by a fuzzy number $\overline{r} \in \overline{\overline{\mathbb{R}}},\ \overline{\overline{\mathbb{R}}} = \overline{\overline{\mathbb{Z}}} \cup \overline{\overline{\mathbb{Q}}} \cup \overline{\overline{\mathbb{Q^{\prime}}}}$, where if $r \in \mathbb{Q}^{\prime}$ or $r \in \mathbb{Q} \backslash \mathbb{Z}$  will be replaced by a triangular fuzzy number because of density of irrational and rational numbers in $\mathbb{R}$ and if $r \in \mathbb{Z}$ will be replaced by a singleton fuzzy set because of non density of integer numbers in $\mathbb{R}$ [1], and then using the representation theorem (resolution principle) to calculate the arithmetic operations act on $\alpha-$cuts of fuzzy numbers [4]. And in 2011, Altai defined the limit fuzzy number of the convergent fuzzy sequence in similar way [2]. Our goal is to establish the theory of fuzzy limits of fuzzy functions depending on Altai's principle, because it is very handy and convenient in the study of the fuzzy arithmetic.

\paragraph*{Representation theorem [4].}Let $A$ be a fuzzy set in $X$ with the membership function $\mu_{A}(x)$. Let $A_{\alpha}$ be the $\alpha-$cuts of $A$ and $\chi_{A_{\alpha}}$ be the characteristic function of the crisp set $A_{\alpha}, \alpha \in (0,1]$. Then
\begin{align*}
\mu_{A}(x) = \sup_{\alpha \in (0,1]} \left( \alpha \wedge \chi_{A_{\alpha}}(x) \right), \ x \in X. 
\end{align*} 

\paragraph*{Resolution principle [4].} Let $A$ be a fuzzy set in $X$ and $\alpha A_{\alpha}, \alpha \in (0,1]$ be a special fuzzy set, whose membership function
\begin{align*}
\mu_{\alpha A_{\alpha}}(x) = \left( \alpha \wedge \chi_{A_{\alpha}}(x) \right), \ x \in X. 
\end{align*} 
Also, let 
\begin{align*}
\Lambda_{A} = \left\lbrace \alpha : \mu_{A}(x) = \alpha \ \mathrm{for \ some} \ x \in X \right\rbrace  
\end{align*} 
be the level set of $A$. Then $A$ can be expressed in the form
\begin{align*}
A = \bigcup_{\alpha \in \Lambda_{A}} \left( \alpha A_{\alpha} \right),   
\end{align*} 
where $\bigcup$ denotes the standard fuzzy union. 

\paragraph*{Remark [4].} The essence of representation theorem of fuzzy sets is that a fuzzy set $A$ in $X$ can be retrieved as a union of its $\alpha A_{\alpha}$ sets, $\alpha \in (0,1]$ and the essence of resolution principle is that a fuzzy set $A$ can be decomposed into fuzzy sets $\alpha A_{\alpha}, \alpha \in (0,1]$ . Thus the representation theorem and the resolution principle are the same coin with two sides as both of them essentially tell that a fuzzy set $A$ in $X$ can always be expressed in terms of its $\alpha-$cuts without explicitly resorting to its membership function $\mu_{A}(x)$.

\paragraph*{Proposition [1]} Let $A$ be a fuzzy number, then $A_{\alpha}$ is a closed, convex and compact subset of $\mathbb{R}$, for all $\alpha \in (0,1]$. 

\section{Two-sided fuzzy limits.} 
In this section, definition of the fuzzy limit of fuzzy functions will be introduced and its properties will be considered.

\paragraph{Theorem 2.1.} Let $\Big( \overline{\overline{X}}, \rho \Big) $ and $\left( \overline{\overline{Y}}, d \right)$ be fuzzy metric spaces. Suppose that $f:\overline{\overline{E}} \subset \overline{\overline{X}} \to \overline{\overline{Y}}$ and $\overline{p}$ is a fuzzy limit point of $\overline{\overline{E}}$. If for all $\alpha \in (0,1]$, the bounds of $\alpha-$cut of $f (\overline{x})$ converge to the bounds of $\alpha-$cut of $\overline{L}$, then $f (\overline{x})$ converges to $\overline{L} \in \overline{\overline{Y}}$ as $\overline{x} \to \overline{p}$.
\paragraph{Proof.}  For all $\alpha \in (0,1]$, let $[f_{1}(x_{1,\alpha},x_{2,\alpha}),f_{2}(x_{1,\alpha},x_{2,\alpha})], [L_{1,\alpha}, L_{2,\alpha}]$ be $\alpha-$cuts of $f (\overline{x})$ and $\overline{L}$ respectively, such that for all $\varepsilon > 0$, there exists $\delta_{1}, \delta_{2} > 0$,
\begin{align*}
0 < \rho_{1}((x_{1,\alpha},x_{2,\alpha}),(p_{1,\alpha},p_{2,\alpha})) < \delta_{1} \Rightarrow  d_{1} \left( f_{1}(x_{1,\alpha},x_{2,\alpha}), L_{i,\alpha} \right) < \varepsilon, \\
0 < \rho_{2}((x_{1,\alpha},x_{2,\alpha}),(p_{1,\alpha},p_{2,\alpha})) < \delta_{2} \Rightarrow  d_{2} \left( f_{2}(x_{1,\alpha},x_{2,\alpha}), L_{i,\alpha} \right) < \varepsilon, \\
\end{align*} 
where 
\begin{align*}
d_{1} \left( f_{1}(x_{1,\alpha},x_{2,\alpha}), L_{i,\alpha} \right) = \min \{ d \left( f_{1}(x_{1,\alpha},x_{2,\alpha}), L_{i,\alpha} \right): i = 1,2 \}, \\
d_{2} \left( f_{1}(x_{1,\alpha},x_{2,\alpha}), L_{i,\alpha} \right) = \max \{ d \left( f_{1}(x_{1,\alpha},x_{2,\alpha}), L_{i,\alpha} \right): i = 1,2 \}, \\
\rho_{1}((x_{1,\alpha},x_{2,\alpha}),(p_{1,\alpha},p_{2,\alpha})) = \min \{ \rho((x_{1,\alpha},x_{2,\alpha}),(p_{1,\alpha},p_{2,\alpha})): i =1,2 \}, \\
\rho_{2}((x_{1,\alpha},x_{2,\alpha}),(p_{1,\alpha},p_{2,\alpha})) = \max \{ \rho((x_{1,\alpha},x_{2,\alpha}),(p_{1,\alpha},p_{2,\alpha})): i =1,2 \}.
\end{align*} 
If $f_{*}(x_{1,\alpha},x_{2,\alpha}) \in [f_{1}(x_{1,\alpha},x_{2,\alpha}),f_{2}(x_{1,\alpha},x_{2,\alpha})]$, by the squeeze theorem for functions that 
\begin{align*}
0 < \rho_{*}((x_{1,\alpha},x_{2,\alpha}),(p_{1,\alpha},p_{2,\alpha})) < \delta_{*} \Rightarrow  d_{*} \left( f_{*}(x_{1,\alpha},x_{2,\alpha}), L_{*,\alpha} \right) < \varepsilon, 
\end{align*} 
where $\delta_{*} = \min \{ \delta_{1}, \delta_{2} \}$ and $ L_{*,\alpha} \in [L_{1,\alpha}, L_{2,\alpha}]$. That is, the $\alpha-$cut $[f_{1}(x_{1,\alpha},x_{2,\alpha}),f_{2}(x_{1,\alpha},x_{2,\alpha})]$ of $f (\overline{x})$ converges to the $\alpha-$cut $[L_{1,\alpha}, L_{2,\alpha}]$ of $\overline{L}$ as the $\alpha-$cut $[x_{1,\alpha},x_{2,\alpha}]$ of $\overline{x}$ approaches the $\alpha-$cut $[p_{1,\alpha},p_{2,\alpha}]$ of $\overline{p}$ for all $\alpha \in (0,1]$. By the resolution principle, we complete the proof. $\square$

\paragraph{Theorem 2.2.} Let $\Big( \overline{\overline{X}}, \rho \Big) $ and $\left( \overline{\overline{Y}}, d \right)$ be fuzzy metric spaces. Suppose that $f:\overline{\overline{E}} \subset \overline{\overline{X}} \to \overline{\overline{Y}}$ and $\overline{p}$ is a fuzzy limit point of $\overline{\overline{E}}$. Then $f (\overline{x})$ converges to $\overline{L} \in \overline{\overline{Y}}$ as $\overline{x} \to \overline{p}$ if and only if for all $\alpha \in (0,1]$, for all $\varepsilon > 0$, there exists $\delta > 0$,
\begin{align}
0 < \left\Vert \big( \rho_{1}((x_{1,\alpha},x_{2,\alpha}),(p_{1,\alpha},p_{2,\alpha})), \rho_{2}((x_{1,\alpha},x_{2,\alpha}),(p_{1,\alpha},p_{2,\alpha})) \big) \right\Vert &< \delta \notag\\
\Rightarrow \left\Vert \left(  d_{1} \left( f(x_{1,\alpha},x_{2,\alpha}), L_{i,\alpha} \right), d_{2} \left( f(x_{1,\alpha},x_{2,\alpha}), L_{i,\alpha} \right) \right)  \right\Vert &< \varepsilon,
\end{align}
\paragraph{Proof.} Let $f (\overline{x})$ converge to $\overline{L} \in \overline{\overline{Y}}$ as $\overline{x} \to \overline{p}$. By theorem 2.1, for all $\alpha \in (0,1]$, for all $\varepsilon > 0$, there exists $\delta_{1}, \delta_{2} > 0$,
\begin{align*}
0 < \rho_{1}((x_{1,\alpha},x_{2,\alpha}),(p_{1,\alpha},p_{2,\alpha})) < \delta_{1} \big/ \sqrt{2} \Rightarrow  d_{1} \left( f_{1}(x_{1,\alpha},x_{2,\alpha}), L_{i,\alpha} \right) < \varepsilon \big/ \sqrt{2}, \\
0 < \rho_{2}((x_{1,\alpha},x_{2,\alpha}),(p_{1,\alpha},p_{2,\alpha})) < \delta_{2} \big/ \sqrt{2} \Rightarrow  d_{2} \left( f_{2}(x_{1,\alpha},x_{2,\alpha}), L_{i,\alpha} \right) < \varepsilon \big/ \sqrt{2}.
\end{align*}
Then 
\begin{align*}
0 < \left\Vert \big( \rho_{1}((x_{1,\alpha},x_{2,\alpha}),(p_{1,\alpha},p_{2,\alpha})), \rho_{2}((x_{1,\alpha},x_{2,\alpha}),(p_{1,\alpha},p_{2,\alpha})) \big) \right\Vert &= \\
\left( \left( \rho_{1}((x_{1,\alpha},x_{2,\alpha}),(p_{1,\alpha},p_{2,\alpha})) \right)^{2} + \left( \rho_{2}((x_{1,\alpha},x_{2,\alpha}),(p_{1,\alpha},p_{2,\alpha})) \right)^{2} \right)^{1/2} &< \delta
\end{align*}
where $\delta = \min \{ \delta_{1}, \delta_{2} \}$, implies
\begin{align*}
\left\Vert \left( d_{1} \left( f_{1}(x_{1,\alpha},x_{2,\alpha}), L_{i,\alpha} \right), d_{2} \left( f_{2}(x_{1,\alpha},x_{2,\alpha}), L_{i,\alpha} \right) \right)  \right\Vert &= \\
\left( \left( d_{1} \left( f_{1}(x_{1,\alpha},x_{2,\alpha}), L_{i,\alpha} \right) \right)^{2} + \left( d_{2} \left( f_{2}(x_{1,\alpha},x_{2,\alpha}), L_{i,\alpha} \right) \right)^{2} \right)^{1/2}  &< \varepsilon.
\end{align*}
Now suppose (2.1) is given. Since
\begin{align*}
\rho_{1}((x_{1,\alpha},x_{2,\alpha}),(p_{1,\alpha},p_{2,\alpha})) &\leq \left\Vert \big( \rho_{1}((x_{1,\alpha},x_{2,\alpha}),(p_{1,\alpha},p_{2,\alpha})), \rho_{2}((x_{1,\alpha},x_{2,\alpha}),(p_{1,\alpha},p_{2,\alpha})) \big) \right\Vert; \\
\rho_{2}((x_{1,\alpha},x_{2,\alpha}),(p_{1,\alpha},p_{2,\alpha})) &\leq \left\Vert \big( \rho_{1}((x_{1,\alpha},x_{2,\alpha}),(p_{1,\alpha},p_{2,\alpha})), \rho_{2}((x_{1,\alpha},x_{2,\alpha}),(p_{1,\alpha},p_{2,\alpha})) \big) \right\Vert
\end{align*}
and
\begin{align*}
d_{1} \left( f_{1}(x_{1,\alpha},x_{2,\alpha}), L_{i,\alpha} \right) \leq \left\Vert \left(  d_{1} \left( f_{1}(x_{1,\alpha},x_{2,\alpha}), L_{i,\alpha} \right), d_{2} \left( f_{2}(x_{1,\alpha},x_{2,\alpha}), L_{i,\alpha} \right) \right)  \right\Vert; \\
d_{2} \left( f_{2}(x_{1,\alpha},x_{2,\alpha}), L_{i,\alpha} \right) \leq \left\Vert \left(  d_{1} \left( f_{1}(x_{1,\alpha},x_{2,\alpha}), L_{i,\alpha} \right), d_{2} \left( f_{2}(x_{1,\alpha},x_{2,\alpha}), L_{i,\alpha} \right) \right)  \right\Vert.
\end{align*}
Then
\begin{align*}
0 < \rho_{1}((x_{1,\alpha},x_{2,\alpha}),(p_{1,\alpha},p_{2,\alpha})) < \delta &\Rightarrow d_{1} \left( f_{1}(x_{1,\alpha},x_{2,\alpha}), L_{i,\alpha} \right) < \varepsilon; \\
0 < \rho_{2}((x_{1,\alpha},x_{2,\alpha}),(p_{1,\alpha},p_{2,\alpha})) < \delta &\Rightarrow d_{2} \left( f_{2}(x_{1,\alpha},x_{2,\alpha}), L_{i,\alpha} \right) < \varepsilon. \ \square
\end{align*}

\paragraph{Remark 2.1.}We will call $\overline{L}$ in theorem 2.2 by the fuzzy limit of $f$ at $\overline{p}$ and write it as
\begin{align}
f(\overline{p}) = \overline{L} = \lim \limits_{\overline{x} \to \overline{p}} f(\overline{x}).
\end{align} 

\paragraph{Examples 2.1.}
\begin{enumerate}
\item To find the limit of $\overline{f}(\overline{x})= \frac{\overline{x}^3 -\overline{4}}{\overline{x}^2+\overline{1}}$, as $\overline{x} \to (0,\frac{1}{2},1)$. We have, by the resolution principle, for all $\alpha \in (0,1]$, the $\alpha-$cut 
\begin{align*}
\frac{ \left[ x_{1,\alpha},x_{2,\alpha} \right]^{3} - [4,4]}{ \left[ x_{1,\alpha},x_{2,\alpha} \right]^{2} + [1,1]} &= \left[ \min_{i,j,k = 1,2} \left\lbrace \frac{x_{i,\alpha}x_{j,\alpha}x_{k,\alpha} - 4}{x_{i,\alpha}x_{j,\alpha}+1} \right\rbrace, \max_{i,j,k = 1,2} \left\lbrace \frac{x_{i,\alpha}x_{j,\alpha}x_{k,\alpha} - 4}{x_{i,\alpha}x_{j,\alpha} + 1} \right\rbrace \right] \mathrm{of} \ \overline{f}(\overline{x})
\end{align*}
has the limit
\begin{align*}
\left[ \lim \limits_{ \substack{x_{1,\alpha} \to \frac{1}{2}\alpha \\ x_{2,\alpha} \to {1 - \frac{1}{2}\alpha}}} \min_{i,j,k = 1,2} \left\lbrace \frac{x_{i,\alpha}x_{j,\alpha}x_{k,\alpha} - 4}{x_{i,\alpha}x_{j,\alpha}+1} \right\rbrace, \lim \limits_{ \substack{x_{1,\alpha} \to \frac{1}{2}\alpha \\ x_{2,\alpha} \to {1 - \frac{1}{2}\alpha}}} \max_{i,j,k = 1,2} \left\lbrace \frac{x_{i,\alpha}x_{j,\alpha}x_{k,\alpha} - 4}{x_{i,\alpha}x_{j,\alpha}+1} \right\rbrace \right] .
\end{align*}
Taking the union of above $\alpha-$cut we get the limit of the function. 
\item If $\overline{f}(\overline{x}) = \overline{x} + \overline{b}, \overline{x} \in \overline{\overline{\mathbb{R}}}$, then $ \lim \limits_{\overline{x} \to \overline{p}} \overline{f}(\overline{x}) = \overline{f}(\overline{p}) $ 
because, by the resolution principle, for all $\alpha \in (0,1]$, for all $\varepsilon > 0$, there exists an $\delta > 0$,
\begin{align*}
0 < \left\Vert \left( \left| x_{1,\alpha} - p_{2,\alpha} \right|, \left| x_{2,\alpha} - p_{1,\alpha} \right| \right) \right\Vert &< \delta  \\
\Rightarrow \left\Vert \left( \left| f_{1}(x_{1,\alpha},x_{2,\alpha}) - f_{2}(p_{1,\alpha}, p_{2,\alpha}) \right|, \left| f_{2}(x_{1,\alpha},x_{2,\alpha}) - f_{1}(p_{1,\alpha}, p_{2,\alpha}) \right| \right) \right\Vert &= \\
\left\Vert \left( \left| \left( x_{1,\alpha} + b_{1,\alpha} \right)  - \left( p_{2,\alpha} + b_{2,\alpha} \right)  \right|, \left| \left( x_{2,\alpha} - p_{2,\alpha} \right)  - \left( p_{1,\alpha} + b_{1,\alpha} \right) \right|_{2} \right) \right\Vert &\leq \\
\left\Vert \left( \left| \left( x_{1,\alpha} - p_{2,\alpha} \right)  \right|, \left| x_{2,\alpha} - p_{1,\alpha} \right| \right) \right\Vert + \left\Vert \left( \left| \left( b_{1,\alpha} - b_{2,\alpha} \right)  \right|, \left| b_{2,\alpha} - b_{1,\alpha} \right| \right) \right\Vert &< \\
\left\{ \begin{array}{ll} \delta &	\mbox{,$\mathrm{if} \  \overline{b} \in \overline{\overline{\mathbb{Z}}},$}  \\
\delta + \left\Vert \left( \left| \left( b_{1,\alpha} - b_{2,\alpha} \right)  \right|, \left| b_{2,\alpha} - b_{1,\alpha} \right| \right) \right\Vert & \mbox{,$ \mathrm{if} \  \overline{b} \not\in \overline{\overline{\mathbb{Z}}} $}. \\
\end{array} 
\right.
\end{align*}
\item If $\overline{f}(\overline{x})= \overline{x}^2+\overline{x}-\overline{3}, \overline{x} \in \overline{\overline{\mathbb{R}}} $, then $ \lim \limits_{\overline{x} \to \overline{1}} \overline{f}(\overline{x}) = - \overline{1}$ because, by the resolution principle, for all $\alpha \in (0,1]$, for all $\varepsilon > 0$, there exists $0 < \delta \leq 1$,
\begin{align*}
0 < \left\Vert \left( \left| x_{1,\alpha} - 1 \right|, \left| x_{2,\alpha} - 1 \right| \right) \right\Vert < \delta  \Rightarrow \left\Vert \left( \left| f_{1}(x_{1,\alpha},x_{2,\alpha}) - f_{2}(1,1) \right|, \left| f_{2}(x_{1,\alpha},x_{2,\alpha}) - f_{1}(1,1) \right| \right) \right\Vert &= \\
\left\Vert \left( \left| y_{1,\alpha} + x_{1,\alpha} - 2 \right|, \left| y_{2,\alpha} + x_{2,\alpha} - 2 \right| \right) \right\Vert &< \sqrt{32} \delta = \varepsilon
\end{align*}
where $y_{1,\alpha} = \min \{ x^{2}_{1,\alpha},x_{1,\alpha}x_{2,\alpha},x^{2}_{2,\alpha} \}; y_{2,\alpha} = \max \{ x^{2}_{1,\alpha},x_{1,\alpha}x_{2,\alpha},x^{2}_{2,\alpha} \}$
and 
\begin{align*} 
\left| y_{1,\alpha} + x_{1,\alpha} - 2 \right| &\leq \left| x_{1,\alpha} - 1 \right| \left| x_{1,\alpha} + 2 \right| < \left( \left| x_{1,\alpha} \right| + 2 \right) \delta < 4\delta, \mathrm{if} y_{1,\alpha} = x^{2}_{1,\alpha} ;
\\
\left| y_{1,\alpha} + x_{1,\alpha} - 2 \right| &\leq \left| x_{1,\alpha} - 1 \right| \left| x_{2,\alpha} + 1 \right| + \left| x_{2,\alpha} - 1 \right| < \left( \left| x_{2,\alpha} \right| + 1 \right) \delta + \delta < 4 \delta, \mathrm{if} y_{1,\alpha} = x_{1,\alpha} x_{2,\alpha};
\\
 \left|y_{1,\alpha} + x_{1,\alpha} - 2 \right| &\leq \left| x^{2}_{2,\alpha} - 1 \right| + \left| x_{1,\alpha} - 1 \right| < \left( \left| x_{2,\alpha} \right| + 1 \right) \delta + \delta < 4\delta, \mathrm{if} y_{1,\alpha} = x^{2}_{2,\alpha};
\\
\left| y_{2,\alpha} + x_{2,\alpha} - 2 \right| &\leq \left| x^{2}_{1,\alpha} - 1 \right| + \left| x_{2,\alpha} - 1 \right| < \left( \left| x_{1,\alpha} \right| + 1 \right) \delta + \delta < 4\delta, \mathrm{if} y_{2,\alpha} = x^{2}_{1,\alpha};
\\
\left| y_{2,\alpha} + x_{2,\alpha} - 2 \right| &\leq \left| x_{2,\alpha} - 1 \right| \left| x_{1,\alpha} + 1 \right| + \left| x_{1,\alpha} - 1 \right| < \left( \left| x_{1,\alpha} \right| + 1 \right) \delta + \delta < 4 \delta, \mathrm{if} y_{2,\alpha} = x_{1,\alpha} x_{2,\alpha};
\\
\left|y_{2,\alpha} + x_{2,\alpha} - 2 \right| &\leq \left| x_{2,\alpha} - 1 \right| \left| x_{2,\alpha} + 2 \right| < \left( \left| x_{2,\alpha} \right| + 2 \right) \delta < 4\delta, \mathrm{if} y_{2,\alpha} = x^{2}_{2,\alpha}.
\end{align*}
Set $\delta = \min \left\lbrace 1,\varepsilon / \sqrt{32} \right\rbrace $, we complete the proof.
\end{enumerate}
Now, we can consider basic properties of fuzzy limits of fuzzy functions and prove them depending on the above theorems.

\paragraph{Theorem 2.3.} The fuzzy limit of a fuzzy function is unique if it exists. 
\paragraph{Proof.} Suppose $f:\overline{\overline{E}} \subset \overline{\overline{X}} \to \overline{\overline{Y}} $ and $\overline{p} \in \overline{\overline{X}}$ is a fuzzy limit point of $\overline{\overline{E}}$. Assume that $\lim \limits_{\overline{x} \to \overline{p}} f(\overline{x}) = \overline{L}; \lim \limits_{\overline{x} \to  \overline{p}} f(\overline{x}) = \overline{M}.$
So, by the resolution principle, for all $\alpha \in (0,1]$, for all $\varepsilon > 0$, there exit $\delta_{1}, \delta_{2} > 0$, such that
\begin{align*}
0 < \left\Vert \left( \rho_{1}((x_{1,\alpha},x_{2,\alpha}),(p_{1,\alpha},p_{2,\alpha})), \rho_{2}((x_{1,\alpha},x_{2,\alpha}),(p_{1,\alpha},p_{2,\alpha})) \right) \right\Vert &< \delta_{1} \\
\Rightarrow \left\Vert \left( d_{1} \left( f_{1}(x_{1,\alpha},x_{2,\alpha}), L_{i,\alpha} \right), d_{2} \left( f_{2}(x_{1,\alpha},x_{2,\alpha}), L_{i,\alpha} \right) \right) \right\Vert &< \frac{\varepsilon}{2}; 
\\
0 < \left\Vert \left( \rho_{1}((x_{1,\alpha},x_{2,\alpha}),(p_{1,\alpha},p_{2,\alpha})), \rho_{2}((x_{1,\alpha},x_{2,\alpha}),(p_{1,\alpha},p_{2,\alpha})) \right) \right\Vert &< \delta_{2} \\
\Rightarrow \left\Vert \left( d_{1} \left( f_{1}(x_{1,\alpha},x_{2,\alpha}), M_{i,\alpha} \right), d_{2} \left( f_{2}(x_{1,\alpha},x_{2,\alpha}), M_{i,\alpha} \right) \right) \right\Vert &< \frac{\varepsilon}{2} .
\end{align*}
Let $\delta = \min \{ \delta_{1}, \delta_{2} \}$. Then, for all $\alpha \in (0,1]$, the $\alpha-$cut $\left[ p_{1,\alpha},p_{2,\alpha} \right]$ of $\overline{p}$ satisfies
\begin{align*}
0 < \left\Vert \left( \rho_{1}((x_{1,\alpha},x_{2,\alpha}),(p_{1,\alpha},p_{2,\alpha})), \rho_{2}((x_{1,\alpha},x_{2,\alpha}),(p_{1,\alpha},p_{2,\alpha})) \right) \right\Vert &< \delta \\
\Rightarrow \left\Vert \left( d_{1} \left( L_{i,\alpha}, M_{i,\alpha} \right), d_{2} \left( L_{i,\alpha}, M_{i,\alpha} \right) \right) \right\Vert \leq \left\Vert \left( d_{1} \left( L_{i,\alpha}, f_{1}(x_{1,\alpha},x_{2,\alpha}) \right), d_{2} \left( L_{i,\alpha}, f_{2}(x_{1,\alpha},x_{2,\alpha}) \right) \right) \right\Vert &+ \\
\left\Vert \left( d_{1} \left( f_{1}(x_{1,\alpha},x_{2,\alpha}), M_{i,\alpha} \right), d_{2} \left( f_{2}(x_{1,\alpha},x_{2,\alpha}), M_{i,\alpha} \right) \right) \right\Vert &< \varepsilon,
\end{align*}
where 
\begin{align*}
d_{1} \left( L_{i,\alpha}, M_{i,\alpha} \right) = \min \{ d \left( L_{i,\alpha}, M_{i,\alpha} \right): i =1,2 \}, d_{2} \left( L_{i,\alpha}, M_{i,\alpha} \right) = \max \{ d \left( L_{i,\alpha}, M_{i,\alpha} \right): i =1,2 \}. \ \square
\end{align*}

\paragraph{Theorem 2.4.} Let $f:\overline{ \overline{E}} \subset \overline{\overline{X}} \to \overline{\overline{Y}}$ and $\overline{p}$ be a fuzzy limit point of $\overline{\overline{E}}$. Then $\lim \limits_{\overline{x} \to \overline{p}} f(\overline{x}) = \overline{L}$ if and only if $\lim \limits_{n \to \infty} f(\overline{p}_{n}) = \overline{L}$ for every fuzzy sequence $\overline{p}_{n}$ in $\overline{ \overline{E}}$ such that $ \overline{p}_{n} \neq \overline{p}, \ \lim \limits_{n \to \infty} \overline{p}_{n} = \overline{p}$.
\paragraph{Proof.} Suppose that $\lim \limits_{\overline{x} \to \overline{p}} f(\overline{x}) = \overline{L}$ holds. By the resolution principle, for all $\alpha \in (0,1]$, for all  $\varepsilon > 0$, there exists $\delta > 0$, 
\begin{align*}
0 < \left\Vert \left( \rho_{1}((x_{1,\alpha},x_{2,\alpha}),(p_{1,\alpha},p_{2,\alpha})), \rho_{2}((x_{1,\alpha},x_{2,\alpha}),(p_{1,\alpha},p_{2,\alpha})) \right) \right\Vert &< \delta \\ 
\Rightarrow \left\Vert \left( d_{1} \left( f_{1}(x_{1,\alpha},x_{2,\alpha}), L_{i,\alpha} \right), d_{2} \left( f_{2}(x_{1,\alpha},x_{2,\alpha}), L_{i,\alpha} \right) \right) \right\Vert &< \varepsilon.
\end{align*}
Since $ \overline{p}_{n} \to \overline{p}$, then for all $\alpha \in (0,1]$, there exits $N \in \mathbb{N}$ such that for $n > N$, 
\begin{align*}
0 < \left\Vert \left( \rho_{1}((p_{n,1,\alpha}, p_{n,2,\alpha}),(p_{1,\alpha},p_{2,\alpha})), \rho_{2}((p_{n,1,\alpha}, p_{n,2,\alpha}),(p_{1,\alpha},p_{2,\alpha})) \right) \right\Vert &< \delta \\
\Rightarrow \left\Vert \left( d_{1} \left( f_{1}(p_{n,1,\alpha}, p_{n,2,\alpha}), L_{i,\alpha} \right), d_{2} \left( f_{2}(p_{n,1,\alpha}, p_{n,2,\alpha}), L_{i,\alpha} \right) \right) \right\Vert &< \varepsilon.
\end{align*} 
Conversely, assume $\lim \limits_{n \to \infty} f(\overline{p}_{n}) = \overline{L}$ but $\lim \limits_{\overline{x} \to \overline{p}} f(\overline{x}) \neq \overline{L}$. That is, there exists $\varepsilon_{o} > 0$, such that for every $\delta > 0$, that 
\begin{align*}
0 < \left\Vert \left( \rho_{1}((x_{1,\alpha},x_{2,\alpha}),(p_{1,\alpha},p_{2,\alpha})), \rho_{2}((x_{1,\alpha},x_{2,\alpha}),(p_{1,\alpha},p_{2,\alpha})) \right) \right\Vert &< \delta \\
 \mathrm{but} \ \left\Vert \left(  d_{1} \left( f_{1}(x_{1,\alpha},x_{2,\alpha}), L_{i,\alpha} \right), d_{2} \left( f_{2}(x_{1,\alpha},x_{2,\alpha}), L_{i,\alpha} \right) \right) \right\Vert &> \varepsilon_{o} 
\end{align*}
Taking $\delta = \frac{1}{n}, n \in \mathbb{N}$, there is a $\overline{p}_{n}$ in $\overline{ \overline{E}}$ such that 
\begin{align*}
0 < \left\Vert \left( \rho_{1}((p_{n,1,\alpha}, p_{n,2,\alpha}),(p_{1,\alpha},p_{2,\alpha})), \rho_{2}((p_{n,1,\alpha}, p_{n,2,\alpha}),(p_{1,\alpha},p_{2,\alpha})) \right) \right\Vert &< \frac{1}{n} \\ 
\mathrm{but} \ \left\Vert \left( d_{1} \left( f_{1}(p_{n,1,\alpha}, p_{n,2,\alpha}), L_{i,\alpha} \right), d_{2} \left( f_{2}(p_{n,1,\alpha}, p_{n,2,\alpha}), L_{i,\alpha} \right) \right) \right\Vert &> \varepsilon_{o} 
\end{align*}
which contradicts the assumption $\lim \limits_{n \to \infty} f(\overline{p}_{n}) = \overline{L}$. $\square$

\paragraph{Theorem 2.5.}
If $f$ and $g$ are fuzzy functions such that $\lim \limits_{\overline{x} \to \overline{p}} g(\overline{x}) = \overline{L}$ and $\lim \limits_{\overline{u} \to \overline{L}} f(\overline{u}) = f \left( \overline{L} \right) $, then $\lim \limits_{\overline{x} \to \overline{p}} f (\overline{g}(\overline{x})) = f \left( \lim \limits_{\overline{x} \to \overline{p}} g(\overline{x}) \right) = f \left( \overline{L} \right).$ 
\paragraph{Proof.} Since $ f(\overline{u}) \to f \left( \overline{L} \right) $ as $\overline{u} \to \overline{L}$, then by the resolution principle, for all $\alpha \in (0,1]$, for all $\varepsilon > 0$, there exits $\delta > 0$, such that
\begin{align*}
0 < \left\Vert \left( \rho_{1}((u_{1,\alpha},u_{2,\alpha}),(L_{1,\alpha},L_{2,\alpha})), \rho_{2}((u_{1,\alpha},u_{2,\alpha}),(L_{1,\alpha},L_{2,\alpha})) \right) \right\Vert &< \delta \\
\Rightarrow \left\Vert \left( d_{1} \left( f_{1}(u_{1,\alpha},u_{2,\alpha}), f_{i}(L_{1,\alpha}, L_{2,\alpha}) \right), d_{2} \left( f_{2}(u_{1,\alpha},u_{2,\alpha}), f_{i}(L_{1,\alpha}, L_{2,\alpha}) \right) \right) \right\Vert &< \varepsilon.
\end{align*}
Since $ g(\overline{x}) \to \overline{L}$ as $\overline{x} \to \overline{p}$, then by the resolution principle, for all $\alpha \in (0,1]$, there exists $\delta^{\prime} > 0$ such that
\begin{align*}
0 < \left\Vert \left( \sigma_{1}((x_{1,\alpha},x_{2,\alpha}),(p_{1,\alpha},p_{2,\alpha})), \sigma_{2}((x_{1,\alpha},x_{2,\alpha}),(p_{1,\alpha},p_{2,\alpha})) \right) \right\Vert &< \delta^{\prime} \\
\Rightarrow \left\Vert \left( \rho_{1} \left( g_{1} (x_{1,\alpha},x_{2,\alpha}), g_{i}(p_{1,\alpha},p_{2,\alpha}) \right) , \rho_{2} \left( g_{2} (x_{1,\alpha},x_{2,\alpha}), g_{i} (p_{1,\alpha},p_{2,\alpha}) \right) \right) \right\Vert &< \delta.
\end{align*}
Letting $u_{1,\alpha} = g_{1}(x_{1,\alpha},x_{2,\alpha}), u_{2,\alpha} = g_{2}(x_{1,\alpha},x_{2,\alpha})$, we obtain 
\begin{align*}
0 < \left\Vert \left( \sigma_{1}((x_{1,\alpha},x_{2,\alpha}),(p_{1,\alpha},p_{2,\alpha})), \sigma_{2}((x_{1,\alpha},x_{2,\alpha}),(p_{1,\alpha},p_{2,\alpha})) \right) \right\Vert &< \delta^{\prime} \Rightarrow \\
\Big \Vert \big( d_{1} \big( f_{1}(g_{1}(x_{1,\alpha},x_{2,\alpha}),g_{2} (x_{1,\alpha},x_{2,\alpha})), f_{i}(L_{1,\alpha}, L_{2,\alpha}) \big), d_{2} \big( f_{2}(g_{1} (x_{1,\alpha},x_{2,\alpha}),g_{2} (x_{1,\alpha},x_{2,\alpha})), f_{i}(L_{1,\alpha}, L_{2,\alpha}) \big) \big) \big \Vert &< \varepsilon. \ \square
\end{align*}
  
\paragraph{Theorem 2.6.} If $\overline{\overline{E}} \subset \overline{\overline{\mathbb{R}}}$ is a fuzzy metric space, $\overline{p}$ is a fuzzy limit point of $\overline{\overline{E}}$, $f$ and $g$ are fuzzy functions on $\overline{\overline{E}}$, and $\lim \limits_{\overline{x} \to \overline{p}} f(\overline{x})$ and $\lim \limits_{\overline{x} \to \overline{p}} g(\overline{x})$ are exist, then 
\begin{enumerate}
\item $ \lim \limits_{\overline{x} \to \overline{p}} (f(\overline{x}) + g(\overline{x}))= \lim \limits_{\overline{x} \to \overline{p}} f(\overline{x}) + \lim \limits_{\overline{x} \to \overline{p}} g(\overline{x})$
\\
\item $\lim \limits_{\overline{x} \to \overline{p}} \left( \overline{A} f \right)(\overline{x}) = \overline{A}  \lim \limits_{\overline{x} \to \overline{p}} f(\overline{x}), \overline{A} \in \overline{\overline{R}}$
\\
\item $\lim \limits_{\overline{x} \to \overline{p}} (f g)(\overline{x})= \lim \limits_{\overline{x} \to \overline{p}} f(\overline{x}) \lim \limits_{\overline{x} \to \overline{p}} g(\overline{x})$
\\ 
\item $\lim \limits_{\overline{x} \to \overline{p}} \left( \frac{f(\overline{x})}{g(\overline{x})} \right) = \frac{\lim \limits_{\overline{x} \to \overline{p}} f(\overline{x})}{\lim \limits_{\overline{x} \to \overline{p}} g(\overline{x})}$ .
\end{enumerate}
\paragraph{Proof.} For $(1)$ and $(2)$, by the resolution principle, we have
\begin{align*}
\lim \limits_{\overline{x} \to \overline{p}} (f(\overline{x}) + g(\overline{x})) &= \bigcup _{\alpha \in (0,1]} \Bigg( \alpha \Bigg[ \lim \limits_{\substack{x_{1,\alpha} \to p_{1,\alpha} \\ x_{2,\alpha} \to p_{2,\alpha}}} (f_{1}(x_{1,\alpha}, x_{2,\alpha}) + g_{1}(x_{1,\alpha}, x_{2,\alpha})), \lim \limits_{\substack{x_{1,\alpha} \to p_{1,\alpha} \\ x_{2,\alpha} \to p_{2,\alpha}}} (f_{2}(x_{1,\alpha}, x_{2,\alpha}) + g_{2}(x_{1,\alpha}, x_{1,\alpha})) \Bigg] \Bigg) \\
&= \bigcup _{\alpha \in (0,1]} \left( \alpha \left[ \lim \limits_{\substack{x_{1,\alpha} \to p_{1,\alpha} \\ x_{2,\alpha} \to p_{2,\alpha}}} f_{1}(x_{1,\alpha},x_{2,\alpha}),\lim \limits_{\substack{x_{1,\alpha} \to p_{1,\alpha} \\ x_{2,\alpha} \to p_{2,\alpha}}} f_{2}(x_{1,\alpha},x_{2,\alpha}) \right] \right) \\
&+ \bigcup _{\alpha \in (0,1]} \left( \alpha \left[ \lim \limits_{\substack{x_{1,\alpha} \to p_{1,\alpha} \\ x_{2,\alpha} \to p_{2,\alpha}}} g_{1}(x_{1,\alpha}, x_{2,\alpha}), \lim \limits_{\substack{x_{1,\alpha} \to p_{1,\alpha} \\ x_{2,\alpha} \to p_{2,\alpha}}} g_{2}(x_{1,\alpha},x_{2,\alpha}) \right] \right) \\
&= \lim \limits_{\overline{x} \to \overline{p}} \overline{f}(\overline{x}) + \lim \limits_{\overline{x} \to \overline{p}} \overline{g}(\overline{x})
\end{align*}
and   
\begin{align*}
\lim \limits_{\overline{x} \to \overline{p}} \left( \overline{A} f \right)(\overline{x}) &= \bigcup _{\alpha \in (0,1]} \left( \alpha \left[ \lim \limits_{\substack{x_{1,\alpha} \to p_{1,\alpha} \\ x_{2,\alpha} \to p_{2,\alpha}}} F_{1}(x_{1,\alpha}, x_{2,\alpha}), \lim \limits_{\substack{x_{1,\alpha} \to p_{1,\alpha} \\ x_{2,\alpha} \to p_{2,\alpha}}} F_{2}(x_{1,\alpha}, x_{2,\alpha}) \right] \right) \\
&= \bigcup _{\alpha \in (0,1]} \left( \alpha \left[ A_{1,\alpha}, A_{2,\alpha} \right] \right) \bigcup _{\alpha \in (0,1]} \left( \alpha \left[ \lim \limits_{\substack{x_{1,\alpha} \to p_{1,\alpha} \\ x_{2,\alpha} \to p_{2,\alpha}}} f_{1}(x_{1,\alpha}, x_{2,\alpha}), \lim \limits_{\substack{x_{1,\alpha} \to p_{1,\alpha} \\ x_{2,\alpha} \to p_{2,\alpha}}} f_{2}(x_{1,\alpha}, x_{2,\alpha})\right] \right) \\
&= \overline{A}  \lim \limits_{\overline{x} \to \overline{p}} f(\overline{x})
\end{align*}
where
\begin{align*}
F_{1}(x_{1,\alpha}, x_{2,\alpha}) = \min \{ A_{1,\alpha} f_{1}(x_{1,\alpha}, x_{2,\alpha}), A_{1,\alpha} f_{2}(x_{1,\alpha}, x_{2,\alpha}), A_{2,\alpha} f_{1}(x_{1,\alpha}, x_{2,\alpha}), A_{2,\alpha} f_{2}(x_{1,\alpha}, x_{2,\alpha}) \},\\
F_{2}(x_{1,\alpha}, x_{2,\alpha}) = \max \{ A_{1,\alpha} f_{1}(x_{1,\alpha}, x_{2,\alpha}), A_{1,\alpha} f_{2}(x_{1,\alpha}, x_{2,\alpha}), A_{2,\alpha} f_{1}(x_{1,\alpha}, x_{2,\alpha}), A_{2,\alpha} f_{2}(x_{1,\alpha}, x_{2,\alpha}) \}.
\end{align*}
\\
To prove $(3)$, let $\lim \limits_{\overline{x} \to \overline{p}} f (\overline{x}) = \overline{L} $ and $\lim \limits_{\overline{x} \to \overline{p}} g (\overline{x}) = \overline{M} $, then $\lim \limits_{\overline{x} \to \overline{p}} \left[ f (\overline{x}) - \overline{L} \right] = \overline{0} $ and $\lim \limits_{\overline{x} \to \overline{p}} \left[ g (\overline{x}) - \overline{M} \right] = \overline{0} $. By the resolution principle, for all $\alpha \in (0,1]$, for all $\varepsilon > 0$, there exists $\delta > 0$, such that 
\begin{align*}
0 < \left\Vert \left( \left|x_{1,\alpha} - p_{2,\alpha} \right|_{1}, \left|x_{2,\alpha} - p_{1,\alpha} \right| \right) \right\Vert < \delta  \Rightarrow \left\Vert \left( \left| f_{1}(x_{1,\alpha},x_{2,\alpha}) - L_{2,\alpha} \right|, \left| f_{2}(x_{1,\alpha},x_{2,\alpha}) - L_{1,\alpha} \right| \right) \right\Vert &< \varepsilon ;
\\
0 < \left\Vert \left( \left|x_{1,\alpha} - p_{2,\alpha} \right|_{1}, \left|x_{2,\alpha} - p_{1,\alpha} \right| \right) \right\Vert < \delta  \Rightarrow \left\Vert \left( \left| g_{1}(x_{1,\alpha},x_{2,\alpha}) - M_{2,\alpha} \right|, \left| g_{2}(x_{1,\alpha},x_{2,\alpha}) - M_{1,\alpha} \right| \right) \right\Vert &< \varepsilon .
\end{align*}
So,
\begin{align*}
\left\Vert \left( \left| (FG)_{1} \right|, \left| (FG)_{2} \right| \right) \right\Vert \leq \left\Vert \left( \left| F_{1} \right|, \left| F_{2}  \right| \right) \right\Vert \left\Vert \left( \left| G_{1} \right|, \left| G_{2} \right| \right) \right\Vert < \varepsilon.
\end{align*}
where
\begin{align*}
F_{1} = f_{1}(x_{1,\alpha},x_{2,\alpha}) - L_{2,\alpha}, F_{2} = f_{2}(x_{1,\alpha},x_{2,\alpha}) - L_{1,\alpha}, \\
G_{1} = g_{1}(x_{1,\alpha},x_{2,\alpha}) - M_{2,\alpha}, G_{2} = g_{2}(x_{1,\alpha},x_{2,\alpha}) - M_{1,\alpha}, \\
(FG)_{1} = \min \{ F_{1} G_{1}, F_{1} G_{2}, F_{2} G_{1}, F_{2} G_{2} \}, \\
(FG)_{2} = \max \{ F_{1} G_{1}, F_{1} G_{2}, F_{2} G_{1}, F_{2} G_{2} \}.
\end{align*}
That is,
\begin{align*}
\lim \limits_{ \substack{x_{1,\alpha} \to p_{1,\alpha} \\ x_{2,\alpha} \to p_{2,\alpha}}} (FG)_{1} = 0,
\lim \limits_{ \substack{x_{1,\alpha} \to p_{2,\alpha} \\ x_{2,\alpha} \to p_{2,\alpha}}} (FG)_{2}  = 0.
\end{align*}
From properties $(1)$ and $(2)$, if $f_{1}(x_{1,\alpha},x_{2,\alpha}) g_{1}(x_{1,\alpha},x_{2,\alpha}) = \min \{ f_{i}(x_{1,\alpha},x_{2,\alpha}) g_{i}(x_{1},x_{2,\alpha}): i = 1,2 \}$ or $f_{1}(x_{1,\alpha},x_{2,\alpha}) g_{1}(x_{1,\alpha},x_{2,\alpha}) = \max \{ f_{i}(x_{1,\alpha},x_{2,\alpha}) g_{i}(x_{1},x_{2,\alpha}): i = 1,2 \}$, then
\begin{align*}
\lim \limits_{ \substack{x_{1,\alpha} \to p_{1,\alpha} \\ x_{2,\alpha} \to p_{2,\alpha}}} f_{1}(x_{1,\alpha},x_{2,\alpha}) g_{1}(x_{1,\alpha},x_{2,\alpha}) &= \lim \limits_{ \substack{x_{1,\alpha} \to p_{1,\alpha} \\ x_{2,\alpha} \to p_{2,\alpha}}} \Big( \left[ f_{1}(x_{1,\alpha},x_{2,\alpha}) - L_{2,\alpha} \right] \left[ g_{1}(x_{1,\alpha},x_{2,\alpha}) - M_{2,\alpha} \right]  \\
&+ L_{2,\alpha} g_{1}(x_{1,\alpha},x_{2,\alpha}) + M_{2,\alpha} f_{1}(x_{1,\alpha},x_{2,\alpha}) - L_{2,\alpha} M_{2,\alpha} \Big) \\
&= 0 + L_{2,\alpha} M_{2,\alpha} + L_{2,\alpha} M_{2,\alpha} - L_{2,\alpha} M_{2,\alpha} = L_{2,\alpha} M_{2,\alpha}.
\end{align*}
If $f_{1}(x_{1,\alpha},x_{2,\alpha}) g_{2}(x_{1,\alpha},x_{2,\alpha}) = \min \{ f_{i}(x_{1,\alpha},x_{2,\alpha}) g_{i}(x_{1},x_{2,\alpha}): i = 1,2 \}$ or $f_{1}(x_{1,\alpha},x_{2,\alpha}) g_{2}(x_{1,\alpha},x_{2,\alpha}) = \max \{ f_{i}(x_{1,\alpha},x_{2,\alpha}) g_{i}(x_{1},x_{2,\alpha}): i = 1,2 \}$, then
\begin{align*}
\lim \limits_{ \substack{x_{1,\alpha} \to p_{1,\alpha} \\ x_{2,\alpha} \to p_{2,\alpha}}} f_{1}(x_{1,\alpha},x_{2,\alpha}) g_{2}(x_{1,\alpha},x_{2,\alpha}) &= \lim \limits_{ \substack{x_{1,\alpha} \to p_{1,\alpha} \\ x_{2,\alpha} \to p_{2,\alpha}}} \Big( \left[ f_{1}(x_{1,\alpha},x_{2,\alpha}) - L_{2,\alpha} \right] \left[ g_{2}(x_{1,\alpha},x_{2,\alpha}) - M_{1,\alpha} \right]  \\
&+ L_{2,\alpha} g_{2}(x_{1,\alpha},x_{2,\alpha}) + M_{1,\alpha} f_{1}(x_{1,\alpha},x_{2,\alpha}) - L_{2,\alpha} M_{1,\alpha} \Big) \\
&= 0 + L_{2,\alpha} M_{1,\alpha} + L_{2,\alpha} M_{1,\alpha} - L_{2,\alpha} M_{1,\alpha} = L_{2,\alpha} M_{1,\alpha}.
\end{align*}
If $f_{2}(x_{1,\alpha},x_{2,\alpha}) g_{1}(x_{1,\alpha},x_{2,\alpha}) = \min \{ f_{i}(x_{1,\alpha},x_{2,\alpha}) g_{i}(x_{1},x_{2,\alpha}): i = 1,2 \}$ or $f_{2}(x_{1,\alpha},x_{2,\alpha}) g_{1}(x_{1,\alpha},x_{2,\alpha}) = \max \{ f_{i}(x_{1,\alpha},x_{2,\alpha}) g_{i}(x_{1},x_{2,\alpha}): i = 1,2 \}$, then
\begin{align*}
\lim \limits_{ \substack{x_{1,\alpha} \to p_{1,\alpha} \\ x_{2,\alpha} \to p_{2,\alpha}}} f_{2}(x_{1,\alpha},x_{2,\alpha}) g_{1}(x_{1,\alpha},x_{2,\alpha}) &= \lim \limits_{ \substack{x_{1,\alpha} \to p_{1,\alpha} \\ x_{2,\alpha} \to p_{2,\alpha}}} \Big( \left[ f_{2}(x_{1,\alpha},x_{2,\alpha}) - L_{1,\alpha} \right] \left[ g_{1}(x_{1,\alpha},x_{2,\alpha}) - M_{2,\alpha} \right]  \\
&+ L_{1,\alpha} g_{1}(x_{1,\alpha},x_{2,\alpha}) + M_{2,\alpha} f_{2}(x_{1,\alpha},x_{2,\alpha}) - L_{1,\alpha} M_{2,\alpha} \Big) \\
&= 0 + L_{1,\alpha} M_{2,\alpha} + L_{1,\alpha} M_{2,\alpha} - L_{1,\alpha} M_{2,\alpha} = L_{1,\alpha} M_{2,\alpha}.
\end{align*}
If $f_{2}(x_{1,\alpha},x_{2,\alpha}) g_{2}(x_{1,\alpha},x_{2,\alpha}) = \min \{ f_{i}(x_{1,\alpha},x_{2,\alpha}) g_{i}(x_{1},x_{2,\alpha}): i = 1,2 \}$ or $f_{2}(x_{1,\alpha},x_{2,\alpha}) g_{2}(x_{1,\alpha},x_{2,\alpha}) = \max \{ f_{i}(x_{1,\alpha},x_{2,\alpha}) g_{i}(x_{1},x_{2,\alpha}): i = 1,2 \}$, then
\begin{align*}
\lim \limits_{ \substack{x_{1,\alpha} \to p_{1,\alpha} \\ x_{2,\alpha} \to p_{2,\alpha}}} f_{2}(x_{1,\alpha},x_{2,\alpha}) g_{2}(x_{1,\alpha},x_{2,\alpha}) &= \lim \limits_{ \substack{x_{1,\alpha} \to p_{1,\alpha} \\ x_{2,\alpha} \to p_{2,\alpha}}} \Big( \left[ f_{2}(x_{1,\alpha},x_{2,\alpha}) - L_{1,\alpha} \right] \left[ g_{2}(x_{1,\alpha},x_{2,\alpha}) - M_{1,\alpha} \right]  \\
&+ L_{2,\alpha} g_{2}(x_{1,\alpha},x_{2,\alpha}) + M_{1,\alpha} f_{2}(x_{1,\alpha},x_{2,\alpha}) - L_{1,\alpha} M_{1,\alpha} \Big) \\
&= 0 + L_{1,\alpha} M_{1,\alpha} + L_{1,\alpha} M_{1,\alpha} - L_{1,\alpha} M_{1,\alpha} = L_{1,\alpha} M_{1,\alpha}.
\end{align*}
Finally, since $\lim \limits_{\overline{x} \to \overline{p}} g(\overline{x}) = \overline{M} $, then by the resolution principle, for all $\alpha \in (0,1]$, for all $\varepsilon > 0$ there exists $\delta_{1} > 0$ such that
\begin{align*}
0 < \left\Vert \left( \left|x_{1,\alpha} - p_{2,\alpha} \right|, \left|x_{2,\alpha} - p_{1,\alpha} \right| \right) \right\Vert < \delta_{1} \Rightarrow \left\Vert \left( \left| g_{1}(x_{1,\alpha},x_{2,\alpha}) - M_{2,\alpha} \right|,  \left| g_{2}(x_{1,\alpha},x_{2,\alpha}) - M_{1,\alpha} \right| \right) \right\Vert < \varepsilon.
\end{align*}
So,
\begin{align*}
0 < \left\Vert \left( \left|x_{1,\alpha} - p_{2,\alpha} \right|, \left|x_{2,\alpha} - p_{1,\alpha} \right| \right) \right\Vert < \delta_{1} \Rightarrow \left\Vert \left( \left| g_{1}(x_{1,\alpha},x_{2,\alpha}) - M_{2,\alpha} \right|,  \left| g_{2}(x_{1,\alpha},x_{2,\alpha}) - M_{1,\alpha} \right| \right) \right\Vert < \frac{\left\Vert \left( \left| M_{1,\alpha} \right|, \left| M_{2,\alpha} \right| \right) \right\Vert}{2}
\end{align*}
which implies that
\begin{align*}
\left\Vert \left( \left| M_{1,\alpha} \right|, \left| M_{2,\alpha} \right| \right) \right\Vert &\leq \left\Vert \left( \left| g_{1}(x_{1,\alpha},x_{2,\alpha}) \right|, \left| g_{2}(x_{1,\alpha},x_{2,\alpha}) \right| \right) \right\Vert + \left\Vert \left( \left| g_{1}(x_{1,\alpha},x_{2,\alpha}) - M_{2,\alpha} \right|,  \left| g_{2}(x_{1,\alpha},x_{2,\alpha}) - M_{1,\alpha} \right| \right) \right\Vert \\
&< \left\Vert \left( \left| g_{1}(x_{1,\alpha},x_{2,\alpha}) - M_{2,\alpha} \right|,  \left| g_{2}(x_{1,\alpha},x_{2,\alpha}) - M_{1,\alpha} \right| \right) \right\Vert + \frac{\left\Vert \left( \left| M_{1,\alpha} \right|, \left| M_{2,\alpha} \right| \right) \right\Vert}{2}
\end{align*}
Rearranging above, we get
\begin{align*}
\frac{1}{\left\Vert \left( \left| g_{1}(x_{1,\alpha},x_{2,\alpha}) \right|, \left| g_{2}(x_{1,\alpha},x_{2,\alpha}) \right|\right) \right\Vert} < \frac{2}{ \left\Vert \left( \left| M_{1,\alpha} \right|, \left| M_{2,\alpha} \right| \right) \right\Vert}.
\end{align*}
Also, there exists $\delta_{2} > 0$ such that
\begin{align*}
0 < \left\Vert \left( \left|x_{1,\alpha} - p_{2,\alpha} \right|, \left|x_{2,\alpha} - p_{1,\alpha} \right| \right) \right\Vert < \delta_{2} \Rightarrow \left\Vert \left( \left| g_{1}(x_{1,\alpha},x_{2,\alpha}) - M_{2,\alpha} \right|,  \left| g_{2}(x_{1,\alpha},x_{2,\alpha}) - M_{1,\alpha} \right| \right) \right\Vert < \frac{\left\Vert \left( \left| M_{1,\alpha} \right|, \left| M_{2,\alpha} \right| \right) \right\Vert^{2} \varepsilon}{2}.
\end{align*}
Set $\delta = \min \{ \delta_{1},\delta_{2} \}$, then  
\begin{align*}
0 < \left\Vert \left( \left|x_{1,\alpha} - p_{2,\alpha} \right|, \left|x_{2,\alpha} - p_{1,\alpha} \right| \right) \right\Vert < \delta \Rightarrow \left| \frac{1}{\left\Vert \left( \left| g_{1}(x_{1,\alpha},x_{2,\alpha}) \right|, \left| g_{2}(x_{1,\alpha},x_{2,\alpha}) \right| \right) \right\Vert} - \frac{1}{\left\Vert \left( \left| M_{1,\alpha} \right|, \left| M_{2,\alpha} \right| \right) \right\Vert} \right| = \\
\frac{\left| \left\Vert \left( \left| M_{1,\alpha} \right|, \left| M_{2,\alpha} \right| \right) \right\Vert - \left\Vert \left( \left| g_{1}(x_{1,\alpha},x_{2,\alpha}) \right|, \left| g_{2}(x_{1,\alpha},x_{2,\alpha}) \right| \right) \right\Vert \right|}{ \left\Vert \left( \left| g_{1}(x_{1,\alpha},x_{2,\alpha}) \right|, \left| g_{2}(x_{1,\alpha},x_{2,\alpha}) \right| \right) \right\Vert \left\Vert \left( \left| M_{1,\alpha} \right|, \left| M_{2,\alpha} \right| \right) \right\Vert} \leq 
\frac{\left\Vert \left(\left| g_{1}(x_{1,\alpha},x_{2,\alpha}) - M_{2,\alpha} \right|,  \left| g_{2}(x_{1,\alpha},x_{2,\alpha}) - M_{1,\alpha} \right|\right) \right\Vert} {\left\Vert \left( \left| g_{1}(x_{1,\alpha},x_{2,\alpha}) \right|, \left| g_{2}(x_{1,\alpha},x_{2,\alpha}) \right| \right) \right\Vert \left\Vert \left( \left| M_{1,\alpha} \right|, \left| M_{2,\alpha} \right| \right) \right\Vert} < \\
\frac{2}{\left\Vert \left(\left| M_{1,\alpha} \right|, \left| M_{2,\alpha} \right| \right) \right\Vert^{2}} \frac{\left\Vert \left( \left| M_{1,\alpha} \right|, \left| M_{2,\alpha} \right| \right) \right\Vert^{2} \varepsilon}{2} = \varepsilon. \ \square
\end{align*}

\paragraph{Theorem 2.7.} Let $\overline{a} \in \overline{\overline{I}} \subset \overline{\overline{R}}$, where $\overline{\overline{I}}$ is an open fuzzy interval. If $f, g$ are fuzzy functions defined on $\overline{\overline{I}} \backslash \overline{a}$ such that $f(\overline{x}) = g(\overline{x}), \overline{x} \in \overline{\overline{I}} \setminus \overline{a} $ and $f(\overline{x}) \to \overline{L}$ as $\overline{x} \to \overline{a}$, then $\lim \limits_{\overline{x} \to \overline{a}} g(\overline{x}) = \lim \limits_{\overline{x} \to \overline{a}} f(\overline{x})$ 
\paragraph{Proof.} Since $f(\overline{x}) \to \overline{L}$ as $\overline{x} \to \overline{a}$, then by the resolution principle, for all $\alpha \in (0,1]$, for all $\varepsilon > 0$, there exists $\delta > 0$ such that
\begin{align*}
0 < \left\Vert \left( \left| x_{1,\alpha} - p_{2,\alpha} \right|, \left| x_{2,\alpha} - p_{1,\alpha} \right| \right) \right\Vert < \delta \Rightarrow \left\Vert \left( \left| f_{1}(x_{1,\alpha},x_{2,\alpha}) - L_{2,\alpha} \right|, \left| f_{2}(x_{1,\alpha},x_{2,\alpha}) - L_{1,\alpha} \right| \right) \right\Vert < \varepsilon; 
\end{align*}
Since $f(\overline{x}) = g(\overline{x}), \overline{x} \in \overline{\overline{I}} \setminus \overline{a}$, then for all $\alpha \in (0,1]$ that $[ f_{1}(x_{1,\alpha},x_{2,\alpha}), f_{2}(x_{1,\alpha},x_{2,\alpha}) ] = [g_{1}(x_{1,\alpha},x_{2,\alpha}), g_{2}(x_{1,\alpha},x_{2,\alpha})]$.
Thus,
\begin{align*}
0 < \left\Vert \left( \left| x_{1,\alpha} - p_{2,\alpha} \right|, \left| x_{2,\alpha} - p_{1,\alpha} \right| \right) \right\Vert < \delta \Rightarrow \left\Vert \left( \left| g_{1}(x_{1,\alpha},x_{2,\alpha}) - L_{2,\alpha} \right|, \left| g_{2}(x_{1,\alpha},x_{2,\alpha}) - L_{1,\alpha} \right| \right) \right\Vert < \varepsilon. \ \square 
\end{align*}

\paragraph{Theorem 2.8.} \medskip \noindent\textbf{Comparison theorem for fuzzy functions.} Suppose $\overline{a} \in \overline{\overline{I}} \subset \overline{\overline{R}}$, where $\overline{\overline{I}}$ is an open fuzzy interval, and $f,g$ are fuzzy functions defined on $\overline{\overline{I}} \backslash \overline{a}$. If $f$ and $g$ have limits as $\overline{x} \to \overline{a}$ and $f(\overline{x}) \leq g(\overline{x})$ for all $\overline{x} \in \overline{\overline{I}} \backslash \overline{a}$, then $\lim \limits_{\overline{x} \to \overline{p}} f(\overline{x}) \leq \lim \limits_{\overline{x} \to \overline{p}} g(\overline{x}) $.
\paragraph{Proof.} Let $\lim \limits_{\overline{x} \to \overline{p}} f(\overline{x}) = \overline{L}$ and $\lim \limits_{\overline{x} \to \overline{p}} g(\overline{x}) = \overline{M}$, and suppose that $ \overline{L} > \overline{M}$. By the resolution principle, for all $\alpha \in (0,1]$, that $\left[ L_{1,\alpha}, L_{2,\alpha} \right] > \left[ M_{1,\alpha}, M_{2,\alpha} \right] $. Let $\varepsilon_{1} > 0, \varepsilon_{2} > 0, \varepsilon_{1} + \varepsilon_{2} = \frac{1}{2}[L_{1,\alpha} - M_{2,\alpha}]$, there exist $\delta_{1} > 0, \delta_{2} > 0$ such that
\begin{align*}
0 < \left\Vert \left( \left| x_{1,\alpha} - p_{2,\alpha} \right|, \left| x_{2,\alpha} - p_{1,\alpha} \right| \right) \right\Vert < \delta_{1} \Rightarrow \left\Vert \left( \left| f_{1}(x_{1,\alpha},x_{2,\alpha}) - L_{2,\alpha} \right|, \left| f_{2}(x_{1,\alpha},x_{2,\alpha}) - L_{1,\alpha} \right| \right) \right\Vert &< \varepsilon_{1};  
\\
0 < \left\Vert \left( \left| x_{1,\alpha} - p_{2,\alpha} \right|, \left| x_{2,\alpha} - p_{1,\alpha} \right| \right) \right\Vert  < \delta_{2} \Rightarrow \left\Vert \left( \left| g_{1}(x_{1,\alpha},x_{2,\alpha}) - M_{2,\alpha} \right|, \left| g_{2}(x_{1,\alpha},x_{2,\alpha}) - M_{1,\alpha} \right| \right) \right\Vert &< \varepsilon_{2}.
\end{align*}
Letting $\delta = \min \{ \delta_{1}, \delta_{2} \}$, we get
\begin{align*}
0 < \left\Vert \left( \left| x_{1,\alpha} - p_{2,\alpha} \right|, \left| x_{2,\alpha} - p_{1,\alpha} \right| \right) \right\Vert  < \delta \Rightarrow \left( f_{1}(x_{1,\alpha},x_{2,\alpha}) - g_{2}(x_{1,\alpha},x_{2,\alpha}), f_{2}(x_{1,\alpha},x_{2,\alpha}) - g_{1}(x_{1,\alpha},x_{2,\alpha}) \right) &= \\
\left( f_{1}(x_{1,\alpha},x_{2,\alpha}) - L_{2,\alpha}, f_{2}(x_{1,\alpha},x_{2,\alpha}) - L_{1,\alpha} \right) + \left( L_{2,\alpha} - M_{1,\alpha}, L_{1,\alpha} - M_{2,\alpha} \right) &+ \\
\left( M_{1,\alpha} - g_{2}(x_{1,\alpha},x_{2,\alpha}), M_{2,\alpha} - g_{1}(x_{1,\alpha},x_{2,\alpha}) \right) &> \\
\left( L_{2,\alpha} - M_{1,\alpha} - \varepsilon_{1} - \varepsilon_{2}, L_{1,\alpha} - M_{2,\alpha} - \varepsilon_{1} - \varepsilon_{2} \right) > (0,0)
\end{align*}
which contradicts the assumption that $\overline{f}(\overline{x}) \leq \overline{g}(\overline{x})$ for all $\overline{x} \in \overline{\overline{I}} \backslash \overline{p}$. $\square$

\paragraph{Theorem 2.9.} \medskip \noindent\textbf{Squeeze theorem for fuzzy functions.} {\sl Suppose $\overline{p} \in \overline{\overline{I}} \subset \overline{\overline{R}}$, where $\overline{\overline{I}}$ is an open fuzzy interval, and $f,g,h$ are fuzzy functions defined on $\overline{\overline{I}} \backslash \overline{p}$. If $f(\overline{x}) \leq h(\overline{x}) \leq g(\overline{x})$ for all $\overline{x} \in \overline{\overline{I}} \backslash \overline{p}$, and $\lim \limits_{\overline{x} \to \overline{p}} f(\overline{x}) = \lim \limits_{\overline{x} \to \overline{p}} g(\overline{x}) = \overline{L}$ then $\lim \limits_{\overline{x} \to \overline{p}} h(\overline{x}) = \overline{L}$. }
\paragraph{Proof.} Since $\lim \limits_{\overline{x} \to \overline{p}} f(\overline{x}) = \lim \limits_{\overline{x} \to \overline{p}} g(\overline{x}) = \overline{L}$, by the resolution principle, for all $\alpha \in (0,1]$, for all $\varepsilon > 0$, there exist $\delta_{1} > 0; \delta_{2} > 0$ such that 
\begin{align*}
0 < \left\Vert \left( \left| x_{1,\alpha} - p_{1,\alpha} \right|, \left| x_{2,\alpha} - p_{2,\alpha} \right| \right) \right\Vert < \delta_{1} \Rightarrow \left\Vert \left( \left| f_{1}(x_{1,\alpha},x_{2,\alpha}) - L_{2,\alpha} \right|, \left| f_{2}(x_{1,\alpha},x_{2,\alpha}) - L_{1,\alpha} \right| \right) \right\Vert &< \varepsilon; 
\\
0 < \left\Vert \left( \left| x_{1,\alpha} - p_{2,\alpha} \right|, \left| x_{2,\alpha} - p_{1,\alpha} \right| \right) \right\Vert < \delta_{2} \Rightarrow \left\Vert \left( \left| g_{1}(x_{1,\alpha},x_{2,\alpha}) - L_{2,\alpha} \right|, \left| g_{2}(x_{1,\alpha},x_{2,\alpha}) - L_{1,\alpha} \right| \right) \right\Vert &< \varepsilon.
\end{align*}
Since $f(\overline{x}) \leq h(\overline{x}) \leq g(\overline{x})$ for all $\overline{x} \in \overline{\overline{I}} \backslash \overline{a}$, then by resolution principle, for all $\alpha \in (0,1]$, there exists $\delta_{3} > 0$ such that
\begin{align*}
0 < \left\Vert \left( \left| x_{1,\alpha} - p_{2,\alpha} \right|, \left| x_{2,\alpha} - p_{1,\alpha} \right| \right) \right\Vert < \delta^{\prime\prime} \Rightarrow \Big( L_{1,\alpha} - \varepsilon \big/ \sqrt{2}, L_{2,\alpha} - \varepsilon \big/ \sqrt{2} \Big) < \Big(f_{1}(x_{1,\alpha},x_{2,\alpha}),f_{2}(x_{1,\alpha},x_{2,\alpha}) \Big) &\leq \\
\Big( h_{1}(x_{1,\alpha},x_{2,\alpha}),h_{2}(x_{1,\alpha},x_{2,\alpha}) \Big) \leq \Big( g_{1}(x_{1,\alpha},x_{2,\alpha}),g_{2}(x_{1,\alpha},x_{2,\alpha}) \Big) < \Big( L_{1,\alpha} + \varepsilon \big/ \sqrt{2}, L_{2,\alpha} + \varepsilon \big/ \sqrt{2} \Big).
\end{align*}
Choosing $\delta = \min \{ \delta_{1},\delta_{2},\delta_{3} \}$ we have
\begin{align*}
0 < \left\Vert \left( \left| x_{1,\alpha} - p_{2,\alpha} \right|_{1}, \left| x_{2,\alpha} - p_{1,\alpha} \right|_{2} \right) \right\Vert &< \delta \Rightarrow \left( \left| h_{1}(x_{1,\alpha},x_{2,\alpha}) - L_{2,\alpha} \right|, \left| h_{2}(x_{1,\alpha},x_{2,\alpha}) - L_{1,\alpha} \right| \right) &< \left( \frac{\varepsilon}{\sqrt{2}}, \frac{\varepsilon}{\sqrt{2}} \right)
\end{align*}
which completes the proof. $\square $

\section{One-sided fuzzy limit.}
We try in this section to establish the concept of the one-side fuzzy limit of fuzzy functions through the following theorem whose proofs are similar to proofs of theorems 2.1 and 2.2 respectively. 
\paragraph{Theorem 3.1.} Let $f: \overline{\overline{I}} \subset \overline{\overline{\mathbb{R}}} \to \overline{\overline{\mathbb{R}}}$ be a fuzzy function defined on some open fuzzy interval $\overline{\overline{I}}$ with left endpoint $\overline{p}$. Then $f (\overline{x})$ converges to $\overline{L}$ as $\overline{x}$ approaches $\overline{p}$ from the right if for all $\alpha \in (0,1]$, the bounds of $\alpha-$cut of $f (\overline{x})$ converge to the bounds of $\alpha-$cut of $\overline{L}$ as the bounds of $\alpha-$cut of $\overline{x}$ approach from the right to the bounds of $\alpha-$cut of $\overline{p}$.

\paragraph{Theorem 3.2.} Let $f: \overline{\overline{I}} \subset \overline{\overline{\mathbb{R}}} \to \overline{\overline{\mathbb{R}}}$ be a fuzzy function defined on some open fuzzy interval $\overline{\overline{I}}$ with right endpoint $\overline{p}$. Then $f (\overline{x})$ converges to $\overline{L}$ as $\overline{x}$ approaches $\overline{p}$ from the left if for all $\alpha \in (0,1]$, the bounds of $\alpha-$cut of $f (\overline{x})$ converge to the bounds of $\alpha-$cut of $\overline{L}$ as the bounds of $\alpha-$cut of $\overline{x}$ approach from the left to the bounds of $\alpha-$cut of $\overline{p}$.

\paragraph{Theorem 3.3.} Let $f: \overline{\overline{I}} \subset \overline{\overline{\mathbb{R}}} \to \overline{\overline{\mathbb{R}}}$ be a fuzzy function defined on some open fuzzy interval $\overline{\overline{I}}$ with left endpoint $\overline{p}$. Then $f (\overline{x})$ converges to $\overline{L}$ as $\overline{x}$ approaches $\overline{p}$ from the right if and only if for all $\alpha \in (0,1]$, for all $\varepsilon > 0$, there exists $\delta_{1},\delta_{2} > 0$, 
\begin{align}
(0,0) < (x_{1,\alpha} - p_{2,\alpha},x_{2,\alpha} - p_{1,\alpha}) < (\delta_{1}, \delta_{2}) \Rightarrow \left\Vert \left( \left| f_{1}(x_{1,\alpha},x_{2,\alpha}) - L_{2,\alpha} \right|, \left| f_{2}(x_{1,\alpha},x_{2,\alpha}) - L_{1,\alpha} \right| \right) \right\Vert &< \varepsilon.
\end{align}

\paragraph{Theorem 3.4.} Let $f: \overline{\overline{I}} \subset \overline{\overline{\mathbb{R}}} \to \overline{\overline{\mathbb{R}}}$ be a fuzzy function defined on some open fuzzy interval $\overline{\overline{I}}$ with right endpoint $\overline{p}$. Then $f (\overline{x})$ converges to $\overline{L}$ as $\overline{x}$ approaches $\overline{p}$ from the left if and only if for all $\alpha \in (0,1]$, for all $\varepsilon > 0$, there exists $\delta_{1}, \delta_{2} > 0$, 
\begin{align}
(- \delta_{1}, - \delta_{2}) < (x_{1,\alpha} - p_{2,\alpha},x_{2,\alpha} - p_{1,\alpha}) < (0,0) \Rightarrow \left\Vert \left( \left| f_{1}(x_{1,\alpha},x_{2,\alpha}) - L_{2,\alpha} \right|, \left| f_{2}(x_{1,\alpha},x_{2,\alpha}) - L_{1,\alpha} \right| \right) \right\Vert &< \varepsilon.
\end{align}

\paragraph{Remark 3.1.} 
\begin{enumerate}
\item We will call $\overline{L}$ in theorem 3.3 by the right-hand fuzzy limit of $f$ at $\overline{p}$ and write it as
\begin{align}
f(\overline{p}^{+}) = \overline{L} = \lim \limits_{\overline{x} \to \overline{p}^{+}} f(\overline{x})
\end{align} 
if by the resolution principle, for all $\alpha \in (0,1]$, for all  $\varepsilon > 0$, there exists $\delta > 0$ such that
3.1 is satisfied.
\item We will call $\overline{L}$ in theorem 3.4 by the left-hand fuzzy limit of $f$ at $\overline{p}$ and write it as
\begin{align}
f(\overline{p}^{-}) = \overline{L} = \lim \limits_{\overline{x} \to \overline{p}^{-}} f(\overline{x})
\end{align} 
if by the resolution principle, for all $\alpha \in (0,1]$, for all  $\varepsilon > 0$, there exists $\delta > 0$ such that
3.2 is satisfied.
\end{enumerate} 

\paragraph{Examples 3.1.}  
\begin{enumerate} 
\item Both $ \lim \limits_{\overline{x} \to (\frac{1}{4},\frac{1}{3},\frac{1}{2})^{+}} \frac{\overline{1}}{\overline{x} - (\frac{1}{4},\frac{1}{3},\frac{1}{2})}$ and $ \lim \limits_{\overline{x} \to (\frac{1}{4},\frac{1}{3},\frac{1}{2})^{-}} \frac{\overline{1}}{\overline{x} - (\frac{1}{4},\frac{1}{3},\frac{1}{2})}$ do not exist, because by the resolution principle, for all $\alpha \in (0,1]$, for all $\varepsilon > 0; \varepsilon^{\prime} > 0$, there exist an $\delta_{1}, \delta_{2} > 0; \delta^{\prime}_{1}, \delta^{\prime}_{2} > 0$ such that
\begin{align*}
\left( 0, 0 \right) < \left(x_{1,\alpha} - \left( - \frac{1}{6} \alpha + \frac{1}{2} \right) ,x_{2,\alpha} - \left( \frac{1}{12} \alpha + \frac{1}{4} \right) \right) &< (\delta_{1}, \delta_{2}) \\
\Rightarrow \left\Vert \left( \left| \frac{1}{x_{2,\alpha} - \left( \frac{1}{12} \alpha + \frac{1}{4} \right)} - \infty \right|, \left| \frac{1}{x_{1,\alpha} - \left( - \frac{1}{6} \alpha + \frac{1}{4} \right)} - \infty \right| \right) \right\Vert &> \left\Vert (1/ \delta_{2}, 1/\delta_{1}) \right\Vert > \varepsilon; 
\\
\left( -\delta^{\prime}_{1}, -\delta^{\prime}_{2} \right) < \left(x_{1,\alpha} - \left( - \frac{1}{6} \alpha + \frac{1}{2} \right) ,x_{2,\alpha} - \left( \frac{1}{12} \alpha + \frac{1}{4} \right) \right) &< (0, 0) \\
\Rightarrow \left\Vert \left( \left| \frac{1}{x_{2,\alpha} - \left( \frac{1}{12} \alpha + \frac{1}{4} \right)} - \infty \right|, \left| \frac{1}{x_{1,\alpha} - \left( - \frac{1}{6} \alpha + \frac{1}{4} \right)} - \infty \right| \right) \right\Vert &> \left\Vert (1/ \delta^{\prime}_{2}, 1/\delta^{\prime}_{1}) \right\Vert > \varepsilon^{\prime}.
\end{align*}
\item $\lim \limits_{\overline{x} \to \overline{0}^{+}} \exp{\left( \overline{1} / \overline{x} \right)}$ does not exist but $\lim \limits_{\overline{x} \to \overline{0}^{-}} \exp{\left( \overline{1} / \overline{x} \right)}$ exists, because by the resolution principle, for all $\alpha \in (0,1]$, for all $\varepsilon > 0$, there exist an $\delta_{1}, \delta_{2} > 0$ such that
\begin{align*}
&(0,0) < (x_{1,\alpha},x_{2,\alpha}) < ( \delta_{1}, \delta_{2} ) \Rightarrow \\
&\left\Vert \left( \left| \exp{\left( 1 / x_{2,\alpha} \right)} - \infty \right|, \left| \exp{\left( 1 / x_{1,\alpha} \right)} - \infty \right| \right) \right\Vert > \left\Vert \left( \left| \exp{\left( 1 / \delta_{2} \right)} - \infty \right|, \left| \exp{\left( 1 / \delta_{1} \right)} - \infty \right| \right) \right\Vert > \varepsilon,
\end{align*}
and for all $\varepsilon^{\prime} > 0$, there exists $\delta^{\prime}_{1}, \delta^{\prime}_{2} > 0$ such that
\begin{align*}
(-\delta^{\prime}_{1}, -\delta^{\prime}_{2}) < (x_{1,\alpha},x_{2,\alpha}) < (0, 0) \Rightarrow \left\Vert \left( \left| \exp{\left( 1 / x_{2,\alpha} \right)} \right|, \left| \exp{\left( 1 / x_{1,\alpha} \right)} \right| \right) \right\Vert < \left\Vert \left( \left| \exp{\left( -1/\delta^{\prime}_{2} \right)} \right|, \left| \exp{\left( -1 /\delta^{\prime}_{1} \right)} \right| \right) \right\Vert < \varepsilon^{\prime}.
\end{align*}
\item The function   
$ f(\overline{x}) = \left\{ \begin{array}{ll}
  \overline{x}^{2} 		&				 \mbox{ ,$\overline{x} < (\frac{1}{6},\frac{1}{5},\frac{1}{4})$}  \\
(\frac{1}{36},\frac{1}{25},\frac{1}{16}) & \mbox{ ,$(\frac{1}{6},\frac{1}{5},\frac{1}{4}) < \overline{x}$} \\
\end{array} 
\right. $ has both $\lim \limits_{\overline{x} \to (\frac{1}{6},\frac{1}{5},\frac{1}{4})^{-}} f(\overline{x})$ and $\lim \limits_{\overline{x} \to (\frac{1}{6},\frac{1}{5},\frac{1}{4})^{+}} f(\overline{x})$ because 
by the resolution principle, for all $\alpha \in (0,1]$, for all $\varepsilon > 0$, there exists $\delta_{1}, \delta_{2} > 0$ such that
\begin{align*}
&(0,0) < \left( x_{1,\alpha} - \left( \frac{-1}{20} \alpha + \frac{1}{4} \right), x_{2,\alpha} - \left( \frac{1}{30} \alpha + \frac{1}{6} \right)\right) < (\delta_{1},\delta_{2}) \Rightarrow \\
&\left\Vert \left| \left( \frac{11}{900} \alpha + \frac{1}{36} \right) - \left( \frac{11}{900} \alpha + \frac{1}{36} \right) \right|, \left| \left( \frac{-9}{400} \alpha + \frac{1}{16} \right) - \left( \frac{-9}{400} \alpha + \frac{1}{16} \right) \right| \right\Vert < \varepsilon
\end{align*} 
and since for all $\alpha \in (0,1]$, for $\delta^{\prime}_{1}, \delta^{\prime}_{2} > 0$ that
\begin{align*} 
(-\delta^{\prime}_{1}, -\delta^{\prime}_{2}) < \left( x_{1,\alpha} - \left( -\frac{1}{20} \alpha + \frac{1}{4} \right), x_{2,\alpha} - \left( \frac{1}{30} \alpha + \frac{1}{6} \right)\right) < (0,0) 
\end{align*}
leads to
\begin{align*}
\left| x^{2}_{1,\alpha} - \left( - \frac{1}{20} \alpha + \frac{1}{4} \right)^{2} \right| &\leq \left| x_{1,\alpha} - \left( - \frac{1}{20} \alpha + \frac{1}{4} \right) \right| \left| x_{1,\alpha} + \left( - \frac{1}{20} \alpha + \frac{1}{4} \right) \right| < \left( \left| x_{1,\alpha} \right| + 2 \left( - \frac{1}{20} \alpha + \frac{1}{4} \right) \right) \delta_{1}  \\
&< 4 \left( - \frac{1}{20} \alpha + \frac{1}{4} \right) \delta_{1},
\end{align*}
\begin{align*}
\left| x^{2}_{1,\alpha} - \left( \frac{1}{30} \alpha + \frac{1}{6} \right)^{2} \right| &\leq \left| x_{1,\alpha} - \left(- \frac{1}{20} \alpha + \frac{1}{4} \right) \right| \left| x_{1,\alpha} + \left( - \frac{1}{20} \alpha + \frac{1}{4} \right) \right| + \left| \left( - \frac{1}{20} \alpha + \frac{1}{4} \right)^{2} - \left( \frac{1}{30} \alpha + \frac{1}{6} \right)^{2} \right| \\
&< 4 \left( - \frac{1}{20} \alpha + \frac{1}{4} \right) \delta_{1} + \left| \left( - \frac{1}{20} \alpha + \frac{1}{4} \right)^{2} - \left( \frac{1}{30} \alpha + \frac{1}{6} \right)^{2} \right| < 5 \left( - \frac{1}{20} \alpha + \frac{1}{4} \right) \delta_{1},
\end{align*}
\begin{align*}
\left| x^{2}_{2,\alpha} - \left( \frac{1}{30} \alpha + \frac{1}{6} \right)^{2} \right| &\leq \left| x_{2,\alpha} - \left( \frac{1}{30} \alpha + \frac{1}{6} \right) \right| \left| x_{2,\alpha} + \left( \frac{1}{30} \alpha + \frac{1}{6} \right) \right| < \left( \left| x_{2,\alpha} \right| + 2 \left( \frac{1}{30} \alpha + \frac{1}{6} \right) \right) \delta_{2}  \\
&< 4 \left( \frac{1}{30} \alpha + \frac{1}{6} \right) \delta_{2}, 
\end{align*}
\begin{align*}
\left| x^{2}_{2,\alpha} - \left( \frac{-1}{20} \alpha + \frac{1}{4} \right)^{2} \right| &\leq \left| x_{2,\alpha} - \left( \frac{1}{30} \alpha + \frac{1}{6} \right) \right| \left| x_{2,\alpha} + \left( \frac{1}{30} \alpha + \frac{1}{6} \right) \right| + \left| \left( \frac{1}{30} \alpha + \frac{1}{6} \right)^{2} - \left( - \frac{1}{20} \alpha + \frac{1}{4} \right)^{2} \right| \\ 
&< 4 \left( \frac{1}{30} \alpha + \frac{1}{6} \right) \delta_{2} + \left| \left( \frac{1}{30} \alpha + \frac{1}{6} \right)^{2} - \left( - \frac{1}{20} \alpha + \frac{1}{4} \right)^{2} \right| < 5 \left( \frac{1}{30} \alpha + \frac{1}{6} \right) \delta_{2},
\end{align*}
\begin{align*}
\left| x_{1,\alpha}x_{2,\alpha} - \left( - \frac{1}{20} \alpha + \frac{1}{4} \right)^{2} \right| &\leq \left[ \left| x_{1,\alpha} - \left( - \frac{1}{20} \alpha + \frac{1}{4} \right) \right| + \left( - \frac{1}{20} \alpha + \frac{1}{4} \right) \right] \left[ \left| x_{2,\alpha} - \left( \frac{1}{30} \alpha + \frac{1}{6} \right) \right| + \left( \frac{1}{30} \alpha + \frac{1}{6} \right) \right] \\
&+ \left( - \frac{1}{20} \alpha + \frac{1}{4} \right)^{2} < \left[ \delta_{1} + \left( - \frac{1}{20} \alpha + \frac{1}{4} \right) \right] \left[ \delta_{2} + \left( \frac{1}{30} \alpha + \frac{1}{6} \right) \right] + \left( - \frac{1}{20} \alpha + \frac{1}{4} \right)^{2}, 
\end{align*}
\begin{align*}
\left| x_{1,\alpha}x_{2,\alpha} - \left( \frac{1}{30} \alpha + \frac{1}{6} \right)^{2} \right| &\leq \left[ \left| x_{1,\alpha} - \left( - \frac{1}{20} \alpha + \frac{1}{4} \right) \right| + \left( - \frac{1}{20} \alpha + \frac{1}{4} \right) \right] \left[ \left| x_{2,\alpha} - \left( \frac{1}{30} \alpha + \frac{1}{6} \right) \right| + \left( \frac{1}{30} \alpha + \frac{1}{6} \right) \right] \\
&+ \left( \frac{1}{30} \alpha + \frac{1}{6} \right)^{2} < \left[ \delta_{1} + \left( - \frac{1}{20} \alpha + \frac{1}{4} \right) \right] \left[ \delta_{2} + \left( \frac{1}{30} \alpha + \frac{1}{6} \right) \right] + \left( \frac{1}{30} \alpha + \frac{1}{6} \right)^{2}.
\end{align*}
then, by considering above various cases, for all $\varepsilon^{\prime}(\delta_{1},\delta_{2}) > 0$, we get
\begin{align*}
 \left\Vert \left| y_{1,\alpha} - \left( \frac{-1}{20} \alpha + \frac{1}{4} \right)^{2} \right|, \left| y_{2,\alpha} - \left( \frac{1}{30} \alpha + \frac{1}{6} \right)^{2} \right| \right\Vert < \varepsilon^{\prime},
\end{align*}
where $y_{1,\alpha} = \min \{x^{2}_{1,\alpha}, x_{1,\alpha}x_{2,\alpha}, x^{2}_{2,\alpha} \}, y_{2,\alpha} = \max \{x^{2}_{1,\alpha}, x_{1,\alpha}x_{2,\alpha}, x^{2}_{2,\alpha} \} $ and 

\end{enumerate}

\paragraph{Theorem 3.5.} 
$f: \overline{\overline{\mathbb{R}}} \to \overline{\overline{\mathbb{R}}}$ be a fuzzy function, then $
\lim \limits_{\overline{x} \to \overline{p}} f(\overline{x}) = \overline{L}$ if and only if $ 
\overline{L} = \lim \limits_{\overline{x} \to \overline{p}^{-}} f(\overline{x}) = \lim \limits_{\overline{x} \to \overline{p}^{+}} f(\overline{x}). $
\paragraph{Proof.} Suppose that $\lim \limits_{\overline{x} \to \overline{p}} f(\overline{x}) = \overline{L}$. By the resolution principle, for all $\alpha \in (0,1]$, for all $\varepsilon >0$, there exists $\delta > 0$ such that 
\begin{align*}
0 < \left\Vert \left( \left| x_{1,\alpha} - p_{2,\alpha} \right|, \left| x_{2,\alpha} - p_{1,\alpha} \right| \right) \right\Vert < \delta \Rightarrow \left\Vert \left( \left| f_{1}(x_{1,\alpha},x_{2,\alpha}) - L_{2,\alpha} \right|, \left| f_{2} (x_{1,\alpha},x_{2,\alpha}) - L_{1,\alpha} \right| \right) \right\Vert < \varepsilon. 
\end{align*}
Since for all $\alpha \in (0,1]$, that
\begin{align*}
(0,0) < (x_{1,\alpha} - p_{2,\alpha},x_{2,\alpha} - p_{1,\alpha}) < ( \delta_{1} , \delta_{2}) \ \mathrm{and} \ (- \delta_{1}, - \delta_{2}) < (x_{1,\alpha} - p_{2,\alpha},x_{2,\alpha} - p_{1,\alpha}) < (0,0) 
\end{align*}
lead to
\begin{align*}
0 < \left\Vert \left( \left| x_{1,\alpha} - p_{2,\alpha} \right|, \left| x_{2,\alpha} - p_{1,\alpha} \right| \right) \right\Vert < \delta,
\end{align*}
then
\begin{align*}
(0,0) < (x_{1,\alpha} - p_{2,\alpha},x_{2,\alpha} - p_{1,\alpha}) < (\delta_{1} / \sqrt{2}, \delta_{2} / \sqrt{2}) \Rightarrow \left\Vert \left( \left| f_{1}(x_{1,\alpha},x_{2,\alpha}) - L_{2,\alpha} \right|, \left| f_{2}(x_{1,\alpha},x_{2,\alpha}) - L_{1,\alpha} \right| \right) \right\Vert &< \varepsilon; \\
(- \delta_{1} / \sqrt{2}, - \delta_{2} / \sqrt{2}) < (x_{1,\alpha} - p_{2,\alpha},x_{2,\alpha} - p_{1,\alpha}) < (0,0) \Rightarrow \left\Vert \left( \left| f_{1}(x_{1,\alpha},x_{2,\alpha}) - L_{2,\alpha} \right|, \left| f_{2}(x_{1,\alpha},x_{2,\alpha}) - L_{1,\alpha} \right| \right) \right\Vert &< \varepsilon.
\end{align*}
\\
Conversely, suppose $\overline{L} = \lim \limits_{\overline{x} \to \overline{p}-} f(\overline{x}) = \lim \limits_{\overline{x} \to \overline{p}+} f(\overline{x})$ holds. By the resolution principle, for all $\alpha \in (0,1]$, for all $\varepsilon >0$, there exists $\delta > 0$ and $\delta^{\prime} > 0$ such that 
\begin{align*}
(0,0) < (x_{1,\alpha} - p_{2,\alpha},x_{2,\alpha} - p_{1,\alpha}) < (\delta_{1} / \sqrt{2}, \delta_{2} / \sqrt{2}) \Rightarrow \left\Vert \left( \left| f_{1}(x_{1,\alpha},x_{2,\alpha}) - L_{2,\alpha} \right|, \left| f_{2}(x_{1,\alpha},x_{2,\alpha}) - L_{1,\alpha} \right| \right) \right\Vert &< \varepsilon, \\
(- \delta_{1} / \sqrt{2}, - \delta_{2} / \sqrt{2}) < (x_{1,\alpha} - p_{2,\alpha},x_{2,\alpha} - p_{1,\alpha}) < (0,0) \Rightarrow \left\Vert \left( \left| f_{1}(x_{1,\alpha},x_{2,\alpha}) - L_{2,\alpha} \right|, \left| f_{2}(x_{1,\alpha},x_{2,\alpha}) - L_{1,\alpha} \right| \right) \right\Vert &< \varepsilon 
\end{align*}
Set $\delta = \min \{ \delta_{1}, \delta_{2} \}$. Then 
\begin{align*}
0 < \left\Vert \left( \left| x_{1,\alpha} - p_{2,\alpha} \right|_{1}, \left| x_{2,\alpha} - p_{1,\alpha} \right|_{2} \right) \right\Vert < \delta \Rightarrow \left\Vert \left( \left| f(x_{1,\alpha},x_{2,\alpha}) - L_{2,\alpha} \right|, \left| f_{2} (x_{1,\alpha},x_{2,\alpha}) - L_{1,\alpha} \right| \right) \right\Vert < \varepsilon. \ \square
\end{align*}

\paragraph{Examples 3.2.}  
\begin{enumerate}
\item The function $f(\overline{x}) = \frac{\left| \sin(\overline{x}) \right|}{{\sin}(\overline{x})}$ has no fuzzy limit at $\overline{0}$ because by the resolution principle, for all $\alpha \in (0,1]$, we have the $\alpha-$cuts 
$\lim \limits_{\substack{x_{1,\alpha} \to 0^{+} \\ x_{2,\alpha} \to 0^{+}}} \min \left\lbrace \frac{\left| \sin(x_{i,\alpha}) \right|}{\sin(x_{i,\alpha})} : i=1,2 \right\rbrace $ and $\lim \limits_{\substack{x_{1,\alpha} \to 0^{+} \\ x_{2,\alpha} \to 0^{+}}} \max \left\lbrace \frac{\left| \sin(x_{i,\alpha}) \right|}{\sin(x_{i,\alpha})} : i=1,2 \right\rbrace $ give positive values and $\lim \limits_{\substack{x_{1,\alpha} \to 0^{-} \\ x_{2,\alpha} \to 0^{-}}} \min \left\lbrace \frac{\left| \sin(x_{i,\alpha}) \right|}{\sin(x_{i,\alpha})} : i = 1,2 \right\rbrace$ and $\lim \limits_{\substack{x_{1,\alpha} \to 0^{-} \\ x_{2,\alpha} \to 0^{-}}} \max \left\lbrace \frac{\left| \sin(x_{i,\alpha}) \right|}{\sin(x_{i,\alpha})} : i = 1,2 \right\rbrace $ give negative values.
\item The function $f(\overline{x}) = \left\{ \begin{array}{ll}
  \overline{2}\overline{x} + \overline{1} 		&  \mbox{ ,$\overline{x} > \overline{1}$}  \\
  \overline{5}									&  \mbox{ ,$\overline{x} = \overline{1}$}  \\
\overline{7}\overline{x}^{2} - \overline{4}     &  \mbox{ ,$\overline{x} < \overline{1}$} \\
\end{array} 
\right. $ has a fuzzy limit at $\overline{x} = \overline{1}$ because by the resolution principle, for all $\alpha \in (0,1]$, we have the $\alpha-$cuts 
$$ \left[ \lim \limits_{\substack{x_{1,\alpha} \to 1^{+} \\ x_{2,\alpha} \to 1^{+}}} \left( 2x_{1,\alpha} + 1 \right), \lim \limits_{\substack{x_{1,\alpha} \to 1^{+} \\ x_{2,\alpha} \to 1^{+}}} \left( 2x_{2,\alpha} + 1 \right) \right] = \left[ 3,3 \right];$$
$$ \left[ \lim \limits_{\substack{x_{1,\alpha} \to 1^{-} \\ x_{2,\alpha} \to 1^{-}}} \left( 7 y_{1,\alpha} - 4 \right), \lim \limits_{\substack{x_{1,\alpha} \to 1^{-} \\ x_{2,\alpha} \to 1^{-}}} \left( 7 y_{2,\alpha} - 4 \right) \right] = \left[ 3,3 \right].$$
where $y_{1,\alpha} = \min \{ x_{i,\alpha} x_{j,\alpha}: i,j = 1,2 \}; y_{2,\alpha} = \max \{ x_{i,\alpha} x_{j,\alpha}: i,j = 1,2 \}$. Thus, $\lim \limits_{\overline{x} \to \overline{1}^{+}} f(\overline{x}) = \overline{3}; \lim \limits_{\overline{x} \to \overline{1}^{-}} f(\overline{x}) = \overline{3}$, and by theorem 3.5, $\lim \limits_{\overline{x} \to \overline{1}} f(\overline{x}) = \overline{3}$.
\end{enumerate}

\section{Fuzzy limit at infinity.} 

\subsection{ Fuzzy limit as $\overline{x} \to \pm \infty$ } Concept of fuzzy limit of fuzzy function at infinity will be given here through the following theorems whose proofs are similar to proofs of theorems 2.1 and 2.2.

\paragraph{Theorem 4.1.1.} Let $f: \overline{\overline{E}} \subset \overline{\overline{\mathbb{R}}} \to \overline{\overline{\mathbb{R}}}$ be a fuzzy function and $\left( \overline{a},\infty \right) \subseteq \overline{\overline{E}}$ for some $\overline{a} \in \overline{\overline{\mathbb{R}}} $. Then $f (\overline{x})$ converges to $\overline{L} \in \overline{\overline{\mathbb{R}}}$ as $\overline{x}$ approaches $\infty$ if for all $\alpha \in (0,1]$, the bounds of $\alpha-$cut of $f (\overline{x})$ converge to the bounds of $\alpha-$cut of $\overline{L}$ as the bounds of $\alpha-$cut of $\overline{x}$ approach  $\infty$.

\paragraph{Theorem 4.1.2.}Let $f: \overline{\overline{E}} \subset \overline{\overline{\mathbb{R}}} \to \overline{\overline{\mathbb{R}}}$ be a fuzzy function and $\left( - \infty, \overline{a} \right) \subseteq \overline{\overline{E}}$ for some $\overline{a} \in \overline{\overline{\mathbb{R}}} $. Then $f (\overline{x})$ converges to $\overline{L} \in \overline{\overline{\mathbb{R}}}$ as $\overline{x}$ approaches $- \infty$ if for all $\alpha \in (0,1]$, the bounds of $\alpha-$cut of $f (\overline{x})$ converge to the bounds of $\alpha-$cut of $\overline{L}$ as the bounds of $\alpha-$cut of $\overline{x}$ approach  $- \infty$.

\paragraph{Theorem 4.1.3.} Let $f: \overline{\overline{E}} \subset \overline{\overline{\mathbb{R}}} \to \overline{\overline{\mathbb{R}}}$ be a fuzzy function and $\left( \overline{a},\infty \right) \subseteq \overline{\overline{E}}$ for some $\overline{a} \in \overline{\overline{\mathbb{R}}} $. Then $f (\overline{x})$ converges to $\overline{L}$ as $\overline{x}$ approaches $\infty$ if and only if for all $\alpha \in (0,1]$, for all $\varepsilon > 0$, there exists $\overline{K}$ such that the $\alpha-$cuts $\left[ K_{1,\alpha},K_{2,\alpha} \right]$ of $\overline{K}$, $\left[ x_{1,\alpha},x_{2,\alpha} \right]$ of $\overline{x}$ and $ \big[ \big| f_{1}(x_{1,\alpha},x_{2,\alpha}) - L_{2,\alpha} \big|, \big| f_{2}(x_{1,\alpha},x_{2,\alpha}) - L_{1,\alpha} \big| \big]$ of $\big| \overline{f}(\overline{x}) - \overline{L} \big|$ satisfy that $K_{1,\alpha} = K_{1,\alpha}(\varepsilon) > a_{1,\alpha}, K_{2,\alpha} = K_{2,\alpha}(\varepsilon) > a_{2,\alpha}$ and
\begin{align}
\left( x_{1,\alpha}, x_{2,\alpha} \right) > \left( K_{1,\alpha}, K_{2,\alpha} \right) \Rightarrow \left\Vert \left( \big| f_{1}(x_{1,\alpha},x_{2,\alpha}) - L_{2,\alpha} \big|, \big| f_{2}(x_{1,\alpha},x_{2,\alpha}) - L_{1,\alpha} \big| \right)  
 \right\Vert < \varepsilon.
\end{align} 

\paragraph{Theorem 4.1.4.} Let $f: \overline{\overline{E}} \subset \overline{\overline{\mathbb{R}}} \to \overline{\overline{\mathbb{R}}}$ be a fuzzy function and $\left( \overline{a},\infty \right) \subseteq \overline{\overline{E}}$ for some $\overline{a} \in \overline{\overline{\mathbb{R}}} $. Then $f (\overline{x})$ converges to $\overline{L}$ as $\overline{x}$ approaches $\infty$ if and only if for all $\alpha \in (0,1]$, for all $\varepsilon > 0$, there exists $\overline{K}$ such that the $\alpha-$cuts $\left[ K_{1,\alpha},K_{2,\alpha} \right]$ of $\overline{K}$, $\left[ x_{1,\alpha},x_{2,\alpha} \right]$ of $\overline{x}$ and $ \big[ \big| f_{1}(x_{1,\alpha},x_{2,\alpha}) - L_{2,\alpha} \big|, \big| f_{2}(x_{1,\alpha},x_{2,\alpha}) - L_{1,\alpha} \big| \big]$ of $\big| \overline{f}(\overline{x}) - \overline{L} \big|$ satisfy that $K_{1,\alpha} = K_{1,\alpha}(\varepsilon) < a_{1,\alpha}, K_{2,\alpha} = K_{2,\alpha}(\varepsilon) < a_{2,\alpha}$ and
\begin{align}
\left( x_{1,\alpha}, x_{2,\alpha} \right) < \left( K_{1,\alpha}, K_{2,\alpha} \right) \Rightarrow \left\Vert \left( \big| f_{1}(x_{1,\alpha},x_{2,\alpha}) - L_{2,\alpha} \big|, \big| f_{2}(x_{1,\alpha},x_{2,\alpha}) - L_{1,\alpha} \big| \right)  
 \right\Vert < \varepsilon.
\end{align} 

\paragraph{Remark 4.1.1.} The convergence in theorem 4.3 will be denoted as
\begin{align}
\lim \limits_{\overline{x} \to \infty} f(\overline{x}) = \overline{L},
\end{align}
and the convergence in theorem 4.4 will be denoted as
\begin{align}
\lim \limits_{\overline{x} \to - \infty} f(\overline{x}) = \overline{L}.
\end{align} 

\paragraph{Examples 4.1.1.} 
\begin{enumerate}
\item $\lim \limits_{\overline{x} \to \overline{\infty}} \frac{\overline{2} \overline{x}^{2} - \overline{1}}{\overline{1} - \overline{x}^{2}} = - \overline{2}$ because by resolution principle, for all $\alpha \in (0,1]$, the $\alpha-$cut 
\begin{align*}
\left[ \frac{[2,2][x_{1,\alpha},x_{2,\alpha}]^{2}-[1,1]}{[1,1] - [x_{1,\alpha},x_{2,\alpha}]^{2}} \right] = \left[ \min_{i,j = 1,2} \left\lbrace  \frac{2x_{i,\alpha}x_{j,\alpha} -1}{ 1 - x_{i,\alpha}x_{j,\alpha}} \right\rbrace, \max_{i,j = 1,2} \left\lbrace  \frac{2x_{i,\alpha}x_{j,\alpha} -1}{ 1 - x_{i,\alpha}x_{j,\alpha}} \right\rbrace \right]  \mathrm{of} \ \frac{\overline{2} \overline{x}^{2} - \overline{1}}{\overline{1} - \overline{x}^{2}}
\end{align*}
has the limit
\begin{align*}
\left[ \lim \limits_{ \substack{x_{1,\alpha} \to \infty \\ x_{2,\alpha} \to \infty}} \min_{i,j = 1,2} \left\lbrace  \frac{2x_{i,\alpha}x_{j,\alpha} -1}{ 1 - x_{i,\alpha}x_{j,\alpha}} \right\rbrace, \lim \limits_{ \substack{x_{1,\alpha} \to \infty \\ x_{2,\alpha} \to \infty}} \max_{i,j = 1,2} \left\lbrace  \frac{2x_{i,\alpha}x_{j,\alpha} -1}{ 1 - x_{i,\alpha}x_{j,\alpha}} \right\rbrace \right] &= \\
\left[ \lim \limits_{ \substack{x_{1,\alpha} \to \infty \\ x_{2,\alpha} \to \infty}} \min_{i,j = 1,2} \left\lbrace  \frac{2 - 1 / x_{i,\alpha}x_{j,\alpha}}{ -1 + 1/x_{i,\alpha}x_{j,\alpha}} \right\rbrace, \lim \limits_{ \substack{x_{1,\alpha} \to \infty \\ x_{2,\alpha} \to \infty}} \max_{i,j = 1,2} \left\lbrace  \frac{2 - 1 / x_{i,\alpha}x_{j,\alpha}}{ -1 + 1 / x_{i,\alpha}x_{j,\alpha}} \right\rbrace \right] &= \left[ - 2, - 2 \right].
\end{align*}
\item $\lim \limits_{\overline{x} \to \overline{\infty}} \frac{\overline{1}}{\overline{x}} = \overline{0} = \lim \limits_{\overline{x} \to - \overline{\infty}} \frac{\overline{1}}{\overline{x}}$ because by resolution principle, for all $\alpha \in (0,1]$, for all $\varepsilon > 0$, there exist $\alpha-$cuts $\left[ K_{1,\alpha},K_{2,\alpha} \right]$ of $\overline{K} > \overline{0}$ such that
\begin{align*}
\left( x_{1,\alpha}, x_{2,\alpha} \right) > \left( K_{1,\alpha}, K_{2,\alpha} \right)  \Rightarrow \left\Vert \left( \left| \frac{1}{x_{2,\alpha}} \right|, \left| \frac{1}{x_{1,\alpha}} \right| \right) \right\Vert < \left\Vert \left( \left| \frac{1}{K_{2,\alpha}} \right|, \left| \frac{1}{K_{1,\alpha}} \right| \right) \right\Vert < \varepsilon
\end{align*}
and for all $ \varepsilon^{\prime} > 0$, there exist $\alpha-$cuts $\left[ K^{\prime}_{1,\alpha},K^{\prime}_{2,\alpha} \right]$  of $\overline{K^{\prime}} > \overline{0}$ such that
\begin{align*}
\left( x_{1,\alpha}, x_{2,\alpha} \right) < \left( -K^{\prime}_{1,\alpha}, -K^{\prime}_{2,\alpha} \right)  \Rightarrow \left\Vert \left( \left| \frac{1}{x_{2,\alpha}} \right|, \left| \frac{1}{x_{1,\alpha}} \right| \right) \right\Vert <\left\Vert \left( \left| \frac{1}{K^{\prime}_{2,\alpha}} \right|, \left| \frac{1}{K^{\prime}_{1,\alpha}} \right| \right) \right\Vert < \varepsilon^{\prime}.
\end{align*}
\end{enumerate}

\subsection{ Infinity Fuzzy limit } Concept of fuzzy limit of fuzzy function at infinity will be given here through the following theorems whose proofs are similar to proofs of theorems 2.1 and 2.2.

\paragraph{Theorem 4.2.1.} Let $f: \overline{\overline{E}} \subset \overline{\overline{\mathbb{R}}} \to \overline{\overline{\mathbb{R}}}$ be a fuzzy function and $\left( \overline{a},\infty \right) \subseteq \overline{\overline{E}}$ for some $\overline{a} \in \overline{\overline{\mathbb{R}}} $. Then $f(\overline{x})$ converges to $\infty$ as $\overline{x}$ approaches $\overline{a}$ if for all $\alpha \in (0,1]$, the bounds of $\alpha-$cut of $f (\overline{x})$ converge to $\infty$ as the bounds of $\alpha-$cut of $\overline{x}$ approach the bounds of $\alpha-$cut of $\overline{a}$.

\paragraph{Theorem 4.2.2.}Let $f: \overline{\overline{E}} \subset \overline{\overline{\mathbb{R}}} \to \overline{\overline{\mathbb{R}}}$ be a fuzzy function and $\left( - \infty, \overline{a} \right) \subseteq \overline{\overline{E}}$ for some $\overline{a} \in \overline{\overline{\mathbb{R}}} $. Then $f(\overline{x})$ converges to $-\infty$ as $\overline{x}$ approaches $\overline{a}$ if for all $\alpha \in (0,1]$, the bounds of $\alpha-$cut of $f (\overline{x})$ converge to $-\infty$ as the bounds of $\alpha-$cut of $\overline{x}$ approach the bounds of $\alpha-$cut of $\overline{a}$.

\paragraph{Remark 4.2.1.} The convergence in theorem 4.2.1 will be denoted as
\begin{align}
\lim \limits_{\overline{x} \to \overline{a}} f(\overline{x}) = \infty,
\end{align}
and the convergence in theorem 4.2.2 will be denoted as
\begin{align}
\lim \limits_{\overline{x} \to \overline{a}} f(\overline{x}) = - \infty.
\end{align}

\paragraph{Examples 4.2.1.}
\begin{enumerate}
\item $\lim \limits_{\overline{x} \to \overline{0}} \frac{1}{\overline{x}^{2}} = \infty$ because by the resolution principle, for all $\alpha \in (0,1]$, there exists $\alpha-$cuts $\left[ M_{1,\alpha},M_{2,\alpha} \right]$  of $\overline{M} \in \overline{\overline{\mathbb{R}}}$ such that
\begin{align*}
0 < \Vert (x_{1,\alpha},x_{2,\alpha}) \Vert < \delta_{1} \Rightarrow f_{1}(x_{1,\alpha},x_{2,\alpha}) > {1}/{\delta^2_{1}} \Rightarrow \delta_{1} = 1/ M_{2,\alpha}, \\
0 < \Vert (x_{1,\alpha},x_{2,\alpha}) \Vert < \delta_{2} \Rightarrow f_{2}(x_{1,\alpha},x_{2,\alpha}) > {1}/{\delta^2_{2}} \Rightarrow \delta_{2} = 1/ M_{1,\alpha} ,
\end{align*} 
where 
\begin{align*}
f_{1}(x_{1,\alpha},x_{2,\alpha}) = \min \left\lbrace {1}/{x^2_{1,\alpha}}, {1}/{x_{1,\alpha}x_{2,\alpha}}, {1}/{x^2_{2,\alpha}} \right\rbrace, f_{2}(x_{1,\alpha},x_{2,\alpha}) = \max \left\lbrace {1}/{x^2_{1,\alpha}}, {1}/{x_{1,\alpha}x_{2,\alpha}}, {1}/{x^2_{2,\alpha}} \right\rbrace.  
\end{align*}
\item $\lim \limits_{\overline{x} \to \overline{1}-} \frac{\overline{x} + \overline{2}}{\overline{2} \overline{x}^{2} - \overline{3} \overline{x} + \overline{1}} = -\infty$ because by the resolution principle, for all $\alpha \in (0,1]$, there exists $\alpha-$cuts $\left[ M_{1,\alpha},M_{2,\alpha} \right]$  of $\overline{M} < \overline{0}$ such that
\begin{align*}
0 < \Vert ( \vert x_{1,\alpha} -1 \vert, \vert x_{2,\alpha} - 1 \vert ) \Vert < \delta_{1} \Rightarrow f_{1}(x_{1,\alpha},x_{2,\alpha}) = \min_{i,j = 1,2} \left\lbrace \frac{x_{i,\alpha} + 2}{2 x_{i,\alpha}x_{j,\alpha} -3x_{i,\alpha} + 1} \right\rbrace < M_{1, \alpha}, \\
0 < \Vert ( \vert x_{1,\alpha} -1 \vert, \vert x_{2,\alpha} - 1 \vert ) \Vert < \delta_{2} \Rightarrow f_{2}(x_{1,\alpha},x_{2,\alpha}) = \max_{i,j = 1,2} \left\lbrace \frac{x_{i,\alpha} + 2}{2 x_{i,\alpha}x_{j,\alpha} -3x_{i,\alpha} + 1} \right\rbrace < M_{2, \alpha},
\end{align*} 
where $2 x_{i,\alpha}x_{j,\alpha} -3x_{i,\alpha} + 1$ is negative and converges to $0$ as $(x_{1,\alpha},x_{2,\alpha})$ approaches to $(1,1)$ from the left. Therefore, choosing $\delta_{i} \in (0,1), i = 1,2$ such that $(1 - \delta_{1},1 - \delta_{1}) < (x_{1,\alpha},x_{2,\alpha}) < (1,1)$ and $(1 - \delta_{2},1 - \delta_{2}) < (x_{1,\alpha},x_{2,\alpha}) < (1,1)$ imply $2/M_{1} < 2 x_{i,\alpha}x_{j,\alpha} -3x_{i,\alpha} + 1$ and $2/M_{2} < 2 x_{i,\alpha}x_{j,\alpha} -3x_{i,\alpha} + 1$ respectively. Since $(0,0) < (x_{1,\alpha},x_{2,\alpha}) < (1,1)$ imply $(2,2) < (x_{1,\alpha} + 2, x_{2,\alpha} + 2) < (3,3)$, we get the result.
\end{enumerate}
 
\section{Conclusion}
Concept of Limit of function can be generalized to fuzzy limit of fuzzy functions. Basic properties that rule the classical concept of limit of function can be also generalized and proved in light of fuzzy logic and fuzzy sets. Future works like fuzzy continuity, fuzzy derivation and fuzzy integration of fuzzy functions and their properties will be considered depending on concept of fuzzy limit of fuzzy function and its basic properties.


\begin{thebibliography}{}
\bibitem .
Al-Tai A. Q. A. On Fuzzy Markov Chains. Babylon University .M.Sc. Thesis. 2008.

\bibitem .
Al-Tai A. Q. A. On the fuzzy metric spaces. European Journal of Scientific Research, vol. 47, no. 2, pp. 214–229, 2010.

\bibitem .
Al-Tai A. Q. A. On the Fuzzy Convergence. Journal of Applied Mathematics. Hindawi Publishing Corporation, ID 147130, doi:10.1155/2011/147130, 2011.

\bibitem .
Burgin M. Theory of fuzzy limits, Fuzzy sets and systems, vol. 115, no. 3, pp. 433-443, 2000.

\bibitem .
Chandra S.,Bector C. Fuzzy Mathematical Programming and Fuzzy Matrix Games, First Edition. Berlin Heidelberg New York. Germany. 2005.

\bibitem .
Esi A. On some new paranormed sequence spaces of fuzzy numbers defined by Orlicz functions and statistical convergence,” Mathematical Modeling and Analysis, vol. 11, no. 4, pp. 379–388, 2006.

\bibitem .
George A. Veeramani P. On some results in fuzzy metric spaces, Fuzzy Sets and Systems 64, 395–399 (1994).

\bibitem .
George A., Veeramani P. On some results of analysis for fuzzy metric spaces, Fuzzy Sets and Systems 90, 365–368, 1997.

\bibitem .
Hosseini S. B., O'Regan D., Saadati R. Some Results on Intuitionistic Fuzzy Spaces. Iranian Journal of Fuzzy Systems, vol. 4, no. 1, pp. 53-64, 2007.

\bibitem .
Kramosil I., Michalek J. Fuzzy metric and statistical metric spaces. Kybernetika, 336-344, 1975.

\bibitem .
Kwon J.S. On statistical and p-Cesaro convergence of fuzzy numbers. Journal of Computational and Applied Mathematics, vol. 7, no. 1, pp. 757–764, 2003.

\bibitem .
Matloka M. Sequences of fuzzy numbers. BUSEFAL, vol. 28, pp. 28–37, 1986.

\bibitem .
N$\acute{a}$d$\acute{a}$ban S. Fuzzy b-metric spaces. International Journal of Computers Communications and Control, 11(2):273-281, 2016.

\bibitem .
Nanda S. On sequences of fuzzy numbers. Fuzzy Sets and Systems, vol. 33, no. 1, pp. 123–126, 1989.

\bibitem .
Naschie E. On the uncertainty of Cantorian geometry and two-slit experiment. Chaos, Solutions and Fractals, 9:517–29, 1998.

\bibitem .
Zadeh L. Fuzzy Sets : Information and Control, 8, 338-353, 1965.

\end{thebibliography}
\end{document}